\documentclass{london}
\usepackage{amssymb,amsmath,mathrsfs,bbm}
\makeatletter
\renewcommand{\@biblabel}[1]{\bf #1.}
\makeatother

\newcommand{\rf}[1][]{\textup{\eqref{#1}}}
\newcommand{\ts}{\textstyle}
\newcommand{\ds}{\displaystyle}

\newcommand{\ul}{\underline}
\newcommand{\Rda}[2][]{\hat{R}^{#1}_{#2}}
\newcommand{\Rdam}[2][]{(\hat{R}^{-1}){}^{#1}_{#2}}
\newcommand{\kri}{\raisebox{.7mm}{\mbox{${\scriptscriptstyle\circ}$}}}
\newcommand{\f}[1][\bo f]{{\raisebox{.25mm}{\mbox{$^\fettc$}}\! #1}}
\DeclareMathOperator{\id}{id}

\DeclareMathOperator{\Mor}{Mor}
\DeclareMathOperator{\im}{im}
\DeclareMathOperator{\rank}{rank}
\DeclareMathOperator{\tr}{tr}

\newcommand{\sgn}{\mathrm{sgn}}
\newcommand{\op}{\mathrm{op}}
\newcommand{\cop}{\mathrm{cop}}
\newcommand{\wtil}{\widetilde}
\newcommand{\half}{\frac{1}{2}}

\newcommand{\mvec}{{\vec{m}}}   
\newcommand{\nvecp}{{\vec{n}^\prime}}    
   
\newcommand{\nvec}{{\vec{n}}}        
\newcommand{\xvec}{{\vec{x}}}
\newcommand{\xvecp}{{\vec{x}\,^\prime}}

\newcommand{\xinvp}{{\overset{\leftarrow}{x}^\prime}}

\newcommand{\ninvp}{{\overset{\leftarrow}{n}^\prime}}  

\newcommand{\minvp}{{\overset{\leftarrow}{m}^\prime}}

\newcommand{\gtil}{\widetilde{g}}
\newcommand{\sss}{\mathfrak{s}}
\newcommand{\qnum}[1][N]{\lbrack \! \lbrack #1 \rbrack \! \rbrack }
\newcommand{\bo}[1]{\mbox{{\boldmath $ #1$}}}
\newcommand{\dd}{\mathrm{d}}

\newcommand{\ot}{\otimes}
\newcommand{\ott}{{\otimes}}
\newcommand{\ota}{{\otimes}_{\!{\scriptscriptstyle\mathcal{A}}}}
\renewcommand{\SS}{\mathcal{S}}
\newcommand{\A}{\mathcal{A}}
\newcommand{\YY}{\mathcal{Y}}
\newcommand{\AD}{\widehat{\A}}
\newcommand{\BB}{\mathcal{B}}
\newcommand{\OO}{\mathcal{O}}
\newcommand{\CC}{\mathcal{C}}
\newcommand{\CCC}{\mathrm{C}}
\newcommand{\SSS}{\mathrm{S}}
\newcommand{\RR}{\mathcal{R}}
\newcommand{\XX}{\mathcal{X}}
\newcommand{\C}{\mathbbm{C}}
\newcommand{\Z}{\mathbbm{Z}}
\newcommand{\N}{\mathbbm{N}}
\newcommand{\NO}{\mathbbm{N}_0}
\newcommand{\slqn}[1][N]{\mathrm{SL}_q(#1)}
\newcommand{\glqn}[1][N]{\mathrm{GL}_q(#1)}

\newcommand{\oqn}[1][N]{\mathrm{O}_q(#1)}
\newcommand{\soqn}[1][N]{\mathrm{SO}_q(#1)}
\newcommand{\spqn}[1][N]{\mathrm{Sp}_q(#1)}
\newcommand{\nn}{\nonumber}
\newcommand{\dett}{{\mathcal{D}}}
\newcommand{\deti}{\mathcal{D}^{-1}}
\newcommand{\lam}{{\lambda}}

\newcommand{\eps}{\epsilon}
\newcommand{\ve} {\varepsilon}
\newcommand{\sig}{\sigma}

\newcommand{\om}{\omega}
\newcommand{\fettc}{{\scriptstyle{\mathrm{c}}}}
\newcommand{\rac}{{\,\triangleleft\,}}
\newcommand{\fettt}{{t}}
\newcommand{\uhr}{\upharpoonright}
\newcommand{\ii}{{\scriptscriptstyle{\mathrm{Inv}}}}
\renewcommand{\ll}{{\scriptscriptstyle{\mathrm{L}}}}
\newcommand{\rr}{{\scriptscriptstyle{\mathrm{R}}}}
\newcommand{\rt}{\mathfrak{r}}
\newcommand{\dl}{{\Delta_\ll}}
\newcommand{\dr}{{\Delta_\rr}}
\newcommand{\dlr}[1][\rho]{{#1_{(-1)}\ot #1_{(0)}\ot #1_{(1)}}}
\newcommand{\dlra}{{(\id\ot\dr)\dl}}
\newcommand{\Gamm}{\varGamma}
\newcommand{\veeg}{{^\vee}\!\varGamma}
\newcommand{\gvee}{\varGamma^\vee}
\newcommand{\Lam}{\varLambda}
\newcommand{\Lamb}[1][\land]{\varLambda^{#1}}

\newcommand{\gp}{\Gamm_+}
\newcommand{\gm}{\Gamm_-}

\newcommand{\gt}{\Gamm_\tau}
\newcommand{\gmt}{\Gamm_{-\tau}}
\newcommand{\gpz}{\Gamm_{+,z}}
\newcommand{\gmz}{\Gamm_{-,z}}

\newcommand{\gtz}{\Gamm_{\tau,z}}
\newcommand{\gl}{{\Gamm_\ll}}
\newcommand{\gr}{{\Gamm_\rr}}
\newcommand{\gi}{{\Gamm_\ii}}
\newcommand{\glm}[2][k]{\Gamm^{ #1}( #2)}
\newcommand{\glmt}[2][k]{\Gamm^{ #1}_\tau(#2)}
\newcommand{\glmmt}[2][k]{\Gamm^{ #1}_{-\tau}(#2)}
\newcommand{\glmp}[2][k]{\Gamm^{ #1}_+( #2)}
\newcommand{\glmm}[2][k]{\Gamm^{ #1}_-( #2)}
\newcommand{\gd}[1][]{\Gamm^{\land #1}}
\newcommand{\gdt}[1][]{\Gamm^{\land #1}_\tau}
\newcommand{\ld}[1][]{\Lam^{\land #1}}
\newcommand{\gdp}[1][]{\Gamm^{\land #1}_+}

\newcommand{\gdm}[1][]{\Gamm^{\land #1}_-}
\newcommand{\gdl}[1][]{\Gamm^{\land #1}_\ll}
\newcommand{\gdr}[1][]{\Gamm^{\land #1}_\rr}
\newcommand{\gdi}[1][]{{\Gamm^{\land {#1}}_\ii}}
\newcommand{\gten}[1][]{\Gamm^{\ot #1}}
\newcommand{\gtent}[1][]{\Gamm^{\ot #1}_\tau}
\newcommand{\gtenmt} [1][] {\Gamm^{\ot #1}_{-\tau}}

\newcommand{\ldual}[1][]{{g_\ll}_{#1}}
\newcommand{\rdual}[1][]{{g_\rr}_{#1}}
\newcommand{\har}{H}
\newcommand{\deR}{H_{\mathrm{de~R}}}

\newcommand{\Ga}[1]{\Gamm _{#1}}
\newcommand{\Gaw}[2][{}]{\Gamm ^{\wedge #1}_{#2}}
\newcommand{\Gat}[2][{}]{\Gamm ^{\otimes #1}_{#2}}
\newcommand{\mt}{,}
\newcommand{\ctr}[3][{}]{\langle {#2}\mt {#3}\rangle _{#1}}
\newcommand{\ctrp}[2]{\ctr[+]{#1}{#2}}
\newcommand{\ctrm}[2]{\ctr[-]{#1}{#2}}
\newcommand{\ctrpm}[2]{\ctr[\pm ]{#1}{#2}}

\newcommand{\kod}[1][]{\partial^{#1}}
\newcommand{\kodpm}{\kod[\pm]}
\newcommand{\Lap}[2][\tau]{L^{#2}_{#1}}
\newcommand{\Lappm}[1][\tau]{\Lap[#1]{\pm }}
\newcommand{\Lapp}[1][\tau]{\Lap[#1]{+}}
\newcommand{\Lapm}[1][\tau]{\Lap[#1]{-}}

\begin{document}
\authorrunninghead{I.\,HECKENBERGER, A.\,SCH\"ULER}
\titlerunninghead{DE RHAM COHOMOLOGY FOR QUANTUM GROUPS}

\title{De Rham Cohomology\\
and Hodge decomposition\\
for  Quantum Groups
\thanks{{\it 
$2000$ Mathematics Subject Classification} 58B30, 17B37, 58A14, 14F40
\\ \indent
Supported
by the Deutsche Forschungsgemeinschaft}
}
\authors{ISTV\'AN HECKENBERGER \normalsize\theoremfont{and} \authorsize
\authorfont
AXEL SCH\"ULER
}

\abstract{Let $\Gamm=\gtz$ be one of the $N^2$-dimensional bicovariant first order
differential calculi for the quantum groups $\glqn$, $\slqn$,  $\soqn$,
or $\spqn$, where $q$ is a transcendental complex number and $z$ is a regular
parameter. It is shown that the de Rham cohomology  of
Woronowicz' external algebra $\gd$ coincides with the de Rham cohomologies of
its left-coinvariant, its right-coinvariant and its (twosided) coinvariant subcomplexes. In the cases $\glqn$ and $\slqn$ the cohomology
ring is isomorphic to the coinvariant external algebra $\gdi$ and to the
vector space of harmonic forms. We prove a
Hodge decomposition theorem in these cases. The main technical tool is the
spectral decomposition of the quantum Laplace-Beltrami operator. 
}

\begin{article}

\section{Introduction}
About ten years ago a general framework for covariant differential
calculi on Hopf algebras was invented by Woronowicz \cite{a-Woro2}.
Since then covariant first order differential
calculi on quantum groups  have been constructed, studied and classified by many
authors, see for
instance \cite{a-CSchWW1, a-AschCast, a-Jurco1,  a-SchSch1, b-KS}.

In classical differential geometry higher order differential forms naturally
appear. The de Rham cohomology of a compact Lie group $G$ characterises
certain topological properties of $G$.
However, there are only very few papers dealing with higher order differential
calculi and de Rham cohomology on quantum groups. 
Maltsiniotis \cite{a-Malt1} constructed a
multiparameter differential graded bialgebra of $\mathrm{GL}(N)$-type having
the classical dimensions of the bigraded components. Tsygan \cite{a-Tsygan1}
studied the linear $\glqn$-differential calculus in detail.
The de Rham cohomologies of the left-covariant $3D$-calculus and of the
bicovariant $4D_\pm$-calculus were  determined by Woronowicz \cite{a-Woro3}
and Grie\ss l \cite{a-Griessl}, respectively.

The present paper deals with the de Rham cohomology and the Hodge
decomposition of the standard
bicovariant differential calculi on the quantum groups of types A, B, C, and
D. We use Woronowicz' construction of the external algebra. Our first main
result
(Theorem\,\ref{t-co}) says that the embeddings of the left-coinvariant, the
right-coinvariant, and the coinvariant (both left- and right-coinvariant)
subcomplex into the whole complex of differential forms are
quasi-isomorphisms, respectively. This means that their de Rham cohomologies
coincide.
Our second main result (Theorem\,\ref{t-hodge}) is a Hodge decomposition
theorem obtained for the quantum groups of type A. 
The main technical tool is the quantum Laplace-Beltrami operator \cite{a-Heck1}
which is constructed
using the dual pairing of two bicovariant differential calculi. Differential
forms vanishing under the action of the quantum Laplace-Beltrami operator are
called harmonic forms.
If the parameter value $z$ of the differential calculus is regular, then the following three spaces
coincide: the de Rham cohomology ring, the algebra of coinvariant forms, and
the vector space of harmonic forms.
For a special class of non-regular parameter values $z$
however,
there exist additional harmonic forms like
$\dett^{k}\rho$, $k\in\Z$, where $\dett$ is the
 quantum determinant and $\rho$ is a coinvariant differential form.

In case of the orthogonal and symplectic quantum groups there exist harmonic forms which are not closed. Therefore we
only have a restricted Hodge decomposition for elements in the image of the quantum Laplace-Beltrami operator.

Our standing assumption is that the deformation parameter $q$ is a
transcendental complex number. On the one hand this ensures that the
coordinate Hopf algebra $\A$ of the quantum group is cosemisimple and that the
theory of corepresentations of $\A$ corresponds to the representation
theory of the underlying  classical Lie group. On the other
hand it guarantees that there are no other harmonic functions except
from polynomials in the quantum determinant.

This paper is organised as follows. In
Section\,\ref{prelim} we collect some basic definitions and  preliminary facts
needed later. The main result about quasi-isomorphisms is
Theorem\,\ref{t-co}. The Hodge decomposition for $\slqn$ and $\glqn$ is given
in 
Theorem\,\ref{t-hodge}. In Section\,\ref{duality} we prove the isomorphy of the
left-dual and the right-dual Hopf bimodules and we add some properties
of the contraction operator. Section\,\ref{laplace} is devoted to the spectral
decomposition of the quantum Laplace-Beltrami operator. Theorem\,\ref{t-co} is
proven therein. Section\,\ref{s-proof} deals with the duality
of differential and codifferential operators. We use the notion of homomorphic
differential calculi due to Pflaum and Schauenburg \cite{a-PflaumSchau} and 
show that $\gpz$ and $\gmz$ are weakly isomorphic. The proof
of Theorem\,\ref{t-hodge} is given in Section\,\ref{s-proof}.

We close this introduction  by  fixing some assumptions and notations that are
used in the sequel.  All vector spaces,
algebras, bialgebras, etc.\ are meant to be $\C$-vector spaces, unital
$\C$-algebras,
$\C$-bialgebras etc. The linear span of a set $\{a_i| i\in K\}$ is
denoted by $\langle a_i| i\in K\rangle$. The symbol
$\A$ always denotes a  Hopf algebra. We write $\A^\ast$ for the dual
vector space of $\A $ and $\A^\circ$ for the dual Hopf algebra. 
All modules, comodules, and bimodules are assumed to be $\A$-modules, 
$\A$-comodules, and $\A$-bimodules if not specified otherwise. 
The comultiplication, the counit, and the
antipode of $\A $ are denoted by
$\Delta$, $\ve$, and  by $S$, respectively. For a  cosemisimple
Hopf algebra $\A$, let $h$ denote the Haar functional on $\A$.
Let $\bo v=(v^i_j)_{i,j\in K}$  be a corepresentation  of
$\A$.
The coalgebra of matrix elements is denoted by
$\CC(\bo v)=\langle v^i_j\,|\,i,j\in K\rangle$. 
As usual $\bo v^\fettc$ denotes the contragredient corepresentation  of
$\bo v$, where $(v^\fettc)^i_j=S(v^j_i)$.  For the space of intertwiners of two
corepresentations $\bo v$ and $\bo w$  the symbol
$\Mor(\bo v, \bo w)$ is used.  We write $\Mor(\bo v)$ for $\Mor(\bo v,\bo v)$.
We use the same notation $\Mor(\bo f,\bo g)$ for
representations $\bo f$ and $\bo g$ of $\A$. If $A$
is a linear mapping, $A^\fettt$ denotes the transpose (dual) mapping  of $A$ and
$\tr A$ the trace of $A$. 
Lower indices of $A$ always refer
to the components of a tensor product where $A$ acts (`leg numbering').
The unit matrix is denoted by $I$. 
Unless it is explicitly stated otherwise,
we use Einstein's  convention to sum over repeated indices in different factors.
Throughout we assume that $N$ is a positive integer with $N\ge 2$ for quantum groups of type A and $N\ge 3$ for quantum groups of types B, C, and D.
The $q$-numbers we are dealing with are
$\qnum_q:=(q^N-q^{-N})/(q-q^{-1})$. 
We use Sweedler's notation $\Delta(a)=\sum a_{(1)}\ot a_{(2)}$,
$\dl(\rho)=\sum \rho_{(-1)}\ot \rho_{(0)}$
and $\dr(\rho)=\sum \rho_{(0)}\ot \rho_{(1)}$
for the coproduct, for left coactions and for right coactions, respectively.
If $\BB$ is an $\A$-bimodule then the  mapping $b\rac a:=Sa_{(1)}\,b\,a_{(2)}$,
$a\in \A$, $b\in\BB$, is called the
{\em right adjoint action} of $\A$ on $\BB$.

\section{Preliminaries}\label{prelim}
In this section we recall some general notions and facts from the theory of
bicovariant differential calculus \cite{a-Woro2}, which are needed later.
More details and proofs of related or unproven statements can be found in
\cite[Chapter 14]{b-KS}.

\subsection*{Hopf bimodules and bicovariant first order differential calculi}
A {\em Hopf bimodule (bicovariant bimodule)}\/ over $\A$ 
is a bimodule $\Gamm$ together with linear mappings 
$\dl\colon\Gamm\to\A\ot\Gamm$ and $\dr\colon\Gamm\to\Gamm\ot\A$ such that
$(\Gamm,\dl,\dr)$ is a bicomodule, 
$\dl(a\omega b)=\Delta(a)\dl(\omega)\Delta(b)$, 
and $\dr(a\omega b)=\Delta(a)\dr(\omega)\Delta(b)$
for $a,\,b\in\A$ and $\omega\in\Gamm$. 
We call the elements of the vector spaces 
$\Gamm_\ll:=\{\om|\dl(\om)=1\ot\om\}$ and 
 $\Gamm_\rr:=\{\om|\dr(\om)=\om\ot1\}$ {\em left-coinvariant} and  {\em
right-coinvariant}, respectively. The elements of $\gi=\gl\cap\gr$
are called {\em coinvariant}.
The dimension of the Hopf bimodule $\Gamm $ is defined to be the dimension
of the vector space $\gl$. We always assume that $\Gamm $ is finite
dimensional.

For $\rho\in\Gamm$ and $f\in\A^\ast$ we define
$f\ast\rho=\rho_{(0)}f(\rho_{(1)})$ and $\rho\ast f=f(\rho_{(-1)})\rho_{(0)}$. In this
way left and right actions of $\A^\ast$ on $\Gamm$ are defined.
If $\{\om_i\}$ is a basis of the vector space $\gl$ then there exist matrices
$\bo v=(v^i_j)$, $v^i_j\in\A$, and $\bo f=(f^i_j)$, $f^i_j\in\A^\circ$, such
that $\bo v$ is a corepresentation of $\A$, $\bo f$ is a representation of
$\A$, $\dr(\om_j)=\om_i\ot v^i_j$, and $\om_i\rac a=f^i_j(a)\om_j$. We briefly
write $\Gamm=(\bo v,\bo f)$ in this situation.

A {\em first order differential calculus}\/ over $\A$ (FODC for short)  is an
$\A$-bimodule $\Gamm$ with  a linear mapping  $\dd\colon\A\to\Gamm$ which satisfies the Leibniz rule $\dd(ab)=\dd a{\cdot} b+a{\cdot}\dd b$ for
$a,\,b\in\A$,  such that  $\Gamm$ is the linear span of elements $a\, \dd b$ with
$a,b\in\A$.
A FODC $\Gamm$  is called {\em left-covariant}\/ if there exists a linear mapping
$\dl\colon\Gamm\to\A\ot\Gamm$ such that
$\dl(a\dd b)=\Delta(a)(\id\ot\dd)\Delta(b)$
for  $a,b\in\A$.
Similarly, $\Gamm$ is called {\em right-covariant}\/ if there exists a linear
mapping $\dr\colon\Gamm\to\Gamm\ot\A$ such that $\dr(a\dd
b)=\Delta(a)(\dd\ot\id)\Delta(b)$ for $a,b\in\A$. The FODC $\Gamm$ is called
{\em bicovariant}\/  if it is both left- and right-covariant.
In this case $(\Gamm,\dl,\dr)$ is a Hopf 
bimodule.

Let $\Gamm$ be a left-covariant FODC.  A central role plays the mapping
$\om\colon\A\to\gl$
defined by $\om(a)=Sa _{(1)}\dd a_{(2)}$. 
The vector space $\RR =\ker \ve \cap \ker \om $ is a right ideal of
$\ker \ve $. It is called the {\em associated right ideal }
to the left-covariant FODC $\Gamm $.
Suppose that $\{\om_i\,|\,i\in K\}$ is a linear basis of $\gl$. Then there
exist linear functionals
 $X_i\in\A^\circ$, $i\in K$, such that
$\om(a)=\sum_{i\in K}X_i(a)\om_i$ for $a\in\A$. The linear space
$\XX=\langle X_i\,|\,i\in K\rangle $ is called the {\em quantum tangent space}
of $\Gamm$. We have $\dd a=\sum_{i\in K}(X_i\ast a)\om_i$ for $a\in \A $.

\subsection*{Exterior Algebras} We briefly recall Woronowicz' construction of
the external algebra to a given Hopf bimodule $\Gamm$. Obviously
$\gten=\bigoplus_{k\ge0}\gten[k]$ is again a Hopf bimodule. Let $\Lam$
be another Hopf bimodule. Then there exists a unique homomorphism  $\sig\colon\Gamm\ota\Lam\to\Lam\ota\Gamm$ of
Hopf bimodules called  the {\em braiding} with 
\begin{equation}\label{e-sigli}
\sig(\alpha\ota\beta)=\beta_{(0)}\ota(\alpha\rac\beta_{(1)}),\quad
\sig^{-1}(\alpha\ota\beta)=\beta\rac(S^{-1}
\alpha_{(1)})\ota\alpha_{(0)}
\end{equation}
for  $\alpha\in\gl$ and $\beta\in\Lam_\ll$, see
\cite[Subsection~13.1.4]{b-KS}. 
Moreover, $\sig$ satisfies the braid equation $\sig_1\sig_2\sig_1=\sig_2\sig_1\sig_2$.
Let $\SSS_k$ be the symmetric group of $k$ elements and let $s_n$ denote the
simple   transposition of $n$ and $n+1$.  For
$\pi\in\SSS_k$ the expression 
$\pi=s_{i_1}\cdots s_{i_r}$ is called reduced if $r$ is minimal. Since $\sig$
satisfies the braid equation, the bimodule automorphism
$\sig_{\pi}\colon\gten[k]\to\gten[k]$, $\sig_\pi=\sig_{i_1}\cdots\sig_{i_r}$,
does not depend on the choice of the reduced expression for $\pi$. Define the
antisymmetriser $A_k$, $k\ge1$, $A_0=\id$, and the endomorphism $B_{i,j}$, $i+j=k$,  of $\gten[k]$ by   
\begin{equation}
A_k=\sum_{\pi\in\SSS_k}\sgn(\pi)\sig_\pi\quad\text{ and }\quad
B_{i,j}=\sum_{\pi^{-1}\in\CCC_{i,j}}\sgn(\pi)\sig_\pi,
\end{equation}
where $\CCC_{i,j}=\{\pi\in\SSS_{i+j}\,|\,\pi(1)<\cdots<\pi(i),
\pi(i+1)<\cdots<\pi(i+j)\}$ are the shuffle permutations.
We have $
B_{1,i}=\id -\sigma _{1}+\sigma _{1}\sigma _{2}-\ldots
+(-1)^i\sigma _{1}\cdots \sigma _{i}$ and $B_{i,1}=\id -\sigma _{i}+\sigma _{i}\sigma _{i-1}-\ldots
+(-1)^i\sigma _{i}\cdots \sigma _{1}$.
 The above constructions  are possible  for any bimodule isomorphism $\sigma $
which satisfies the braid relation. Therefore, replacing $\sigma $ everywhere
by $\sigma ^{-1}$ the above works as well.
In what follows we will use both kinds of operators and write
$A_k^{\pm}$ and $B_{i,j}^{\pm}$ whenever we are
dealing with $\sigma ^{\pm 1}$.
Now we define {\em Woronowicz' external algebra}
$\gd=\bigoplus_{k\ge0}\gd[k]$ by  $\gd[k]:=\gten[k]{/}\ker A_k$.
It can be shown that $\gd$ is both
a Hopf bimodule and a graded super Hopf
algebra, cf. \cite[Subsection~13.2.2]{b-KS}. 

\subsection*{$\sig$-metrics and contractions}
We recall from  \cite[Section~2]{a-Heck1} the important notion of a 
$\sig$-metric. Let  $(\gp,\gm)$ be a pair of
Hopf bimodules and $\overline{\Gamm }:=\gp\ota\gm\oplus\gm\ota\gp$. A linear
mapping $g\colon\overline{\Gamm }\to\A$ is called a
{\em\mbox{$\sig$-metric}} of
$(\gp,\gm)$ if $g$ is a homomorphism of $\A$-bimodules, $g$ is
non-degenerate, $g$ is $\sig$-symmetric, i.\,e.\ $g{\kri}\sig=g$, and
\begin{equation}\label{e-sig+-}
g_{12}\sig_{23}\sig_{12}=g_{23},\quad g_{23}\sig_{12}\sig_{23}=g_{12}
\end{equation} in
$\Lam\ota \overline{\Gamm }$ and $\overline{\Gamm }\ota\Lam$,
respectively. Here $\Lam$ denotes either
$\gp$ or $\gm$.  The \mbox{$\sig$-metric} $g$ is
called {\em bicovariant} if
\begin{align}
\Delta{\kri} g&=(\id\ot g)\dl=(g\ot\id)\dr
\end{align}
in $\Gamm$. In this paper we are only concerned with bicovariant $\sig$-metrics. 
It is easily seen that \rf[e-sig+-] follows from the bicovariance
\cite[Section~II.\,4]{a-BCDRV}. By abuse of notation we sometimes skip the tensor
sign and write
$g(\xi a,\zeta)=g(\xi, a\zeta)=g(\xi a\ota \zeta)$. We recursively extend
$g$ to an endomomorphism $\gtil$ of the 
Hopf bimodule
$\gtent\ota\gtenmt$, $\tau\in\{+,-\}$, by 
\begin{equation}\label{e-gschl}
\begin{split}
\gtil(\xi,a)&=\xi a,\quad \gtil(a,\zeta)=a\zeta,
\\
\gtil(\xi\ota\rho,\eta\ota\zeta)&=\gtil(\xi\,g(\rho,\eta),\zeta)
\end{split}
\end{equation}
for $\xi\in\gtent$, $\rho\in\gt$, $\zeta\in\gtenmt$, $\eta\in\gmt$, and
$a\in\A$.
If $\xi$ and $\zeta$ are of degree $n$ and $k$, $n\ge k$, then
$\gtil(\xi,\zeta)\in\Gamm_\tau^{\ot n-k}$.
\\
Next we define \textit{contractions} $\ctrpm{\cdot }{\cdot }:
\Gat[k]{\tau }\ota \Gat[l]{-\tau }\to \Gat[|k-l|]{\tau '}$,
$\tau \in \{+,-\}$, where $\tau '=\tau$ for $k\geq l$, and $\tau '=-\tau $
for $k<l$,
by
\begin{equation}\label{eq-ctr}
\begin{aligned}
\ctrpm{\xi }{\zeta}&:=\gtil(B^\pm _{k-l,l}\xi \,\mt\, A^\pm _l\zeta)\quad
\text{for}\quad  k\geq l,\\
\ctrpm{\xi }{\zeta}&:=\gtil(A^\pm _k\xi\, \mt\, B^\pm _{k,l-k}\zeta)\quad
\text{for}\quad k<l.
\end{aligned}
\end{equation}
Since $g$ is, the map $\ctrpm{\cdot}{\cdot}$ is a  homomorphism of Hopf bimodules as well. 
If both $k$ and $l$ are less than two the contraction does
not depend on the sign $\pm$ and we sometimes omit it, $\ctrp{\xi}{\zeta}=\ctrm{\xi}{\zeta}=:\ctr{\xi}{\zeta}$.
\\
The next property shows that the antisymmetriser is symmetric with respect to
$\gtil$. For nonnegative
integers $i,j,k,l$, with $1\leq i+j\leq k,l$, we have 
\begin{align*}
\gtil\kri \bigl( (\id ^{\otimes k-i-j}\ot A^\pm _i \ot \id ^{\otimes j})
\mt \id ^{\otimes l}\bigr) &=
\gtil\kri \bigl( \id ^{\otimes k}\mt (\id ^{\otimes j}\ot A^\pm _i
\ot \id ^{\otimes l-i-j})\bigr) .
\end{align*}
Therefore the definition of $\ctrpm{\cdot }{\cdot }$ can be extended to a
contraction map of {\em exterior algebras} namely to  
$\ctrpm{\cdot }{\cdot }:
\Gaw[k]{\tau }\ota \Gaw[l]{-\tau }\to \Gaw[|k-l|]{\tau '}$,
$\tau \in \{+,-\}$, where $\tau '=\tau$ for $k\geq l$ and $\tau '=-\tau $
for $k<l$.

\subsection*{Differential Calculi on Quantum Groups}
Let $\A$ be the coordinate Hopf algebra $\OO(G_q)$ of one of the quantum groups $\glqn$, $\slqn$,
$\oqn$, $\soqn$, or  $\spqn$ as defined in  \cite[Section~1]{a-FRT}. 
Let $\bo u=(u^i_j)_{i,j=1,\dots,N}$ be the fundamental matrix
corepresentation of $\A$. The corresponding $\Rda{}$-matrix for the A-series
is given in \cite[Subsection~1.2]{a-FRT}. The $\Rda{}$-matrix for the B, C, and D
series as well as the defining the antipode matrix $C$ are given in
\cite[Subsection~1.4]{a-FRT}. Now we define  the invertible diagonal matrix
$D=(d_i\delta_{ij})$, $D\in\Mor(\bo u^{\fettc\fettc},\bo u)$. In case of $\glqn$
and $\slqn$ we set $d_i=q^{N+1-2i}$ and  $\rt=q^N$. For the B, C, and D series set
$D=C^\fettt C^{-1}$ and $\rt=\eps q^{N-\eps}$, where $\eps=1$ in the
orthogonal case and $\eps=-1$ in the symplectic case. Let $\sss=\tr D=\tr
D^{-1}$. For the quantum groups  $\glqn$ and $\oqn$ there exists a nontrivial group-like
central element $\dett$ of $\A$, the quantum determinant. It corresponds to the Young
diagram $(1^N)$ and can be constructed using the $q$-antisymmetric
tensor \cite[Section~5]{a-Hayashi}. Note that $\dett^2=1$ for $\oqn$.

A complex number $x\in\C$ is called {\em admissible} for $\A$ if $x$ is
nonzero for $\glqn$, $x^N=q$ in case $\slqn$, $x^2=1$ in cases $\oqn$,
$\soqn[2n]$, and $\spqn$, and finally,  $x=1$ in case $\soqn[2n+1]$. 
Recall that $\A$ is a  coquasitriangular Hopf algebra (see
\cite[Subsection~10.1]{b-KS}) with universal $r$-form $\bo r_x$,
%
given by $\bo r_x(u^i_j,u^k_l)=x^{-1}\Rda[ki]{jl}$,
where $x$ is an admissible parameter.
The matrices $\bo\ell^+$ and $\bo\ell^-$ of representative functionals on $\A$ are defined by
${\ell^+}^i_j(a)=\bo r_x(a,u^i_j)$ and ${\ell^-}^i_j(a)=\bo
r_y(S(u^i_j),a)$, $x,y$ admissible for $\A$. 
A complex number $z$ is called {\em $2$-admissible} for $\A$ if $z=xy$ for two
admissible numbers $x$ and $y$ for $\A$. 
Throughout the paper we assume $z$ to
be 2-admissible with fixed admissible numbers $x,y$ and  $z=xy$.
Then the Hopf bimodule $\gtz$, $\tau\in\{+,-\}$, is  given as follows. Let $\{\om^\tau_{ij}\,|\,i,j=1,\dots,N\}$  be a linear basis of the space of left-coinvariant forms of
$\gtz$. Define the right coaction $\dr (\rho )$ and the right
adjoint action $\rho\rac a$ for  $\rho\in(\gtz)_\ll$ and $a\in\A$ by 
\begin{align}
\dr(\om^+_{ij})=\om^+_{kl}\ot
(\bo u\ot\bo u^\fettc)^{kl}_{ij},\qquad\dr(\om^-_{ij})=\om^-_{kl}\ot
(\bo u^{\fettc\fettc}\ot\bo u^\fettc)^{kl}_{ij},\label{e-rightcoaction}
\end{align}
\begin{alignat*}{2}
\om^+_{ij}\rac a&=S({\ell^-}^k_i){\ell^+}^j_l(a)\om^+_{kl}&&=\bo r_y(u^k_i,a_{(1)})\bo r_x(a_{(2)},u^j_l)\om^+_{kl},
\\
\om^-_{ij}\rac a&=S^{-1}({\ell^+}^k_i){\ell^-}^j_l(a)\om^-_{kl}&&=\bo r_y(a_{(1)},S(u^k_i))
\bo r_x(S(u^j_l),a_{(2)})\om^-_{kl}.
\end{alignat*}
In shorthand notation we set $\gpz=(\bo u\ot \bo u^\fettc,\bo \ell^{-\fettc}\ot\bo
\ell^+)$ and
$\gmz=(\bo u^{\fettc\fettc}\ot \bo u^\fettc, \f[\bo \ell^+]\ot \bo
\ell^-)$, where $(\f)^i_j:=S^{-1}(f^j_i)$. 
There are unique up to scalars coinvariant 1-forms
$\om^+_0=\sum_{i=1}^N\om^+_{ii}\in\gpz$ and
$\om^-_0=\sum_{i=1}^Nd_i^{-1}\om^-_{ii}{\in\gmz}$.
In particular we have the following right adjoint actions
\begin{xalignat*}{2}
\om^+_{ij}\rac u^m_n&=z^{-1}\Rda[mk]{iv}\Rda[jv]{nl}\,\om^+_{kl},& 
\om^+_0\rac u^m_n&=z^{-1}(\Rda[2]{})^{mk}_{nl}\,\om^+_{kl},
\\
\om^-_{ij}\rac u^m_n&= zd_id_k^{-1}\Rdam[mk]{iv}\Rdam[jv]{nl}\,
\om^-_{kl},&
\om^-_0\rac u^m_n&=zd_k^{-1}(\Rda[-2]{})^{mk}_{nl}\,\om^-_{kl}.
\end{xalignat*}
Recall that  $\Gaw{\tau,z}$ is an   {\em inner} bicovariant differential
calculus with coinvariant 1-form $\om^\tau_0$, that is, the differential
$\dd_\tau$ is given by
\begin{align}\label{eq-diff}
\dd _\tau \rho &=\om^\tau_0 \wedge \rho - (-1)^k\rho \wedge \om^\tau_0,
\quad \rho \in \Gaw[k]{\tau}.
\end{align}
Projecting this equation for $\rho=a\in\A$ to the left-coinvariant part of $\gtz$
we get
\begin{align}\label{e-oma}
\om^\tau(a)=\om^\tau_0\rac a-\ve(a)\om_0^\tau,
\end{align}
where $\om^\tau$ denotes the $\om$-mapping for $\gtz$.  
Let $\{X^\tau_{ij}\,|\,i,j=1,\dots, N\}$ be the  basis of the quantum tangent
space $\XX ^\tau $ of $\gtz$ dual to $\{\om^\tau_{ij}\}$. Explicitly, we have
\begin{xalignat}{2}\label{e-X}
X^+_{ij}&=S({\ell^-}^i_k){\ell^+}^k_j-\delta_{ij},&
X^+_0&=\sss^{-1}(D^{-1})^i_jX^+_{ij},
\\
X^-_{ij}&=d_i^{-1}\bigl(S({\ell^+}^{i}_k){\ell^-}^k_j-\delta_{ij}\bigr),&
X^-_0&=\sss^{-1}{\ts\sum }_i X^-_{ii}.
\end{xalignat}
Here $X^\tau_0$ denotes the dual basis element to $\om^\tau_0$ with respect to the
decomposition of $\dr$ on $\gtz$ into irreducible corepresentations.
The corresponding projections in
\mbox{$\Mor(\bo u\ot\bo u^\fettc)$} and $\Mor(\bo u^{\fettc\fettc}\ot\bo u^\fettc)$
are 
\begin{align}\label{e-P0}
(P^+_0)^{ij}_{kl}=\frac{1}{\sss}d_k^{-1}\delta_{ij}\delta_{kl}\quad\text{ and }\quad
(P^-_0)^{ij}_{kl}=\frac{1}{\sss}d_i^{-1}\delta_{ij}\delta_{kl},
\end{align}
respectively.   Note that $X^\tau_0$ is central in $\A^\circ$ and
$S(X^+_0)=X^-_0$. We set $X_0:=X_0^++X_0^-$.
\\
It was shown in  \cite[Propositions~3.1, 3.3, and 3.4]{a-Heck1}
that the settings
\begin{align}\label{e-metric}
g(a\om^+_{ij},\om^-_{kl})=a D^j_k(D^{-1}){}^l_i\quad\text{ and }\quad 
g(a\om^-_{ij},\om^+_{kl})=a\delta_{jk}\delta_{il}
\end{align}
define a bicovariant $\sig$-metric of $(\gpz,\gmz)$.
Note that
\begin{align}\label{e-ctrom}
g(\om_0^\tau,\om_0^{-\tau})=\sss.
\end{align}
To simplify notations we sometimes write $\gt$ instead of  $\gtz$.

\subsection*{Quantum Laplace-Beltrami Operators}
Our main technical tool to reduce the de Rham cohomology from $\gd$ to the
essentially smaller complex $\gdi$ is the quantum Laplace-Beltrami operator
which is defined below.
For a slightly different notion see also \cite[Section~6]{a-Heck1}.

The mappings $\kodpm_\tau:\Gaw[k]{\tau}\to \Gaw[k-1]{\tau}$, $k\geq 0$, defined
by $\kodpm_\tau(a)=0$ for $a\in \A $ and
\begin{gather}\nn 
\kodpm_\tau\rho =\ctrpm{\rho }{\om^{-\tau}_0}
 +(-1)^k\ctrpm{\om^{-\tau}_0}{\rho }
\end{gather}
for $\rho \in \Gaw[k]{\tau}$, $k>0$, are called
\textit{codifferential operators on} $\Gaw[k]{\tau}$.
The linear mappings $\Lappm :\Gaw[k]{\tau }\to \Gaw[k]{\tau }$
defined by
\begin{equation}\label{e-laplace}
\Lappm :=-\dd_\tau\,\kodpm_\tau+ \kodpm_\tau\,\dd_\tau
\end{equation}
are called {\em quantum Laplace-Beltrami operators}.
The elements of the vector space
\begin{align}\label{e-harmforms}
\har ^\pm (\Gaw[k]{\tau }):=\{\rho \in \Gaw[k]{\tau }\,|\,\Lappm \rho =0\}
\end{align}
are called {\em harmonic $k$-forms}.
If no confusion can arise we sometimes write $\Lappm[]$ instead of $\Lappm $.

On elements $a\in \A$ the Laplace-Beltrami operators defined in \cite{a-Heck1}
and the operators $\Lappm$ coincide, since $\kodpm_\tau(a)=0$.
Therefore we have
\begin{align*}
\Lapp a=\Lapm a&=-2\sss a+\ctr{\om^+_0a}{\om^-_0}
+\ctr{\om^-_0 a}{\om^+_0}
\end{align*}
for $a\in \A $.

\section{Main Results}

Define the sets 
$P_+=\{\lam=(\lam_1,\dots,\lam_N)\,|\,\lam_1\ge\lam_2\ge\dots\ge\lam_N,\,
\lam_i\in\Z\}$ and $P_{++}=\{\lam\in P_+\,|\,\lam_N\ge 0\}$. In case of
$\glqn$ the set $P_+$ parametrises exactly the irreducible corepresentations
of $\A$, see  \cite[Theorem~11.51]{b-KS}, where $(1,0,\dots,0)$, $(0,\dots,0,-1)$, and $(1,1,\dots,1)$
correspond to  $\bo u$, $\bo u^\fettc$, and the determinant $\dett$, respectively.  We extend the notation \cite[Chapter~1]{b-Macdonald} for partitions $\lam\in
P_{++}$  to Young frames with ``negative" columns: for
$\lam\in P_+$  define $|\lam|=\lam_1+\cdots+\lam_N$. For  $i\in\{1,\dots,N\}$
and $j\in\Z$
we write $(i,j)\in\lam$ if $1\le j\le \lam_i$ or $\lam_i<j\le 0$. In this
situation with $x:=(i,j)$ set  $\sgn(x):=1$ if $j\ge1$ and $\sgn(x):=-1$
otherwise. Define the content  $c(x):=j-i$
and  $c(\lam):=\sum_{x\in\lam}\sgn(x)\, c(x)$.  In
particular $\sum_{x\in\lam}\sgn(x)=|\lam|$.

In case of $\glqn$ a complex parameter $z$ is said to be {\em regular} if
for all $\lam,\,\mu\in P_{+}$ the number
\begin{equation} \label{e-hlam}
\begin{split}
F_{\lam\mu}&=(z^{-|\mu|}-2+z^{|\lam|})\qnum_q+(q-q^{-1})\times
\\
&\quad \times\bigl(z^{-|\mu|}\sum_{x\in\mu}
\sgn(x)q^{N+2c(x)}-z^{|\lam|}\sum_{x\in\lam}\sgn(x)q^{-N-2c(x)}\bigr)
\end{split}
\end{equation} 
is nonzero, except for the case $\lam=\mu=(0)$. 
\\
Let $\Lamb=\bigoplus_{k\ge0}\Lamb[k]$ be a differential graded algebra. As
usual
\begin{align}
\deR (\Lamb)&=\bigoplus_{k\ge 0}\deR ^k(\Lamb),&
\deR ^k(\Lamb)&=\ker \dd_k{/}\im \dd_{k-1},
\end{align}
denotes the de Rham cohomology of $\Lamb $.
By the Leibniz rule $\ker \dd=\bigoplus_{k\ge0}\ker\dd_k$
is a subalgebra of $\Lamb$ and $\im \dd=\bigoplus_{k\ge0}\im \dd_k$ is an ideal in
$\ker \dd$. Hence $\deR (\Lamb)$ is an algebra. Since the differential $\dd$ is
bicovariant, the algebras $\gdl$, $\gdr$, and $\gdi$,
and the vector spaces
$\glm[\wedge]{\bo 1,\dett}:=\{\rho\in\gdl\,|\,\dr\rho=\rho\ot\dett\}$ and
$\glm[\wedge]{\dett,\bo 1}:=\{\rho\in\gdr\,|\,\dl\rho=\dett\ot \rho\}$ are
$\Z $-graded differential complexes.

\begin{theorem}\label{t-co}
Suppose that $q$ is a transcendental complex number.\\
{\em (a)}
Let $G_q$ denote one of the quantum groups $\glqn$, $\slqn$,
$\soqn$, or $\spqn$, and $\A=\OO(G_q)$ its coordinate Hopf algebra.
Let $\Gamm$ be one of
the $N^2$-dimensional bicovariant first order differential calculi $\gtz$
over $\A $, where $z$ is $2$-admissible. In the case $\glqn$ we assume in
addition that $z$ is regular.  Then we have canonical isomorphisms
$$
\deR (\gd)\cong \deR (\gdl)\cong \deR (\gdr)\cong \deR (\gdi).
$$
{\em (b)} Let $G_q$ be one of the quantum groups $\glqn$ or $\oqn[2n+1]$,
$\A=\OO(G_q)$, and $\Gamm$ as above.
In the case $\glqn$ we assume that $z^Nq^{-2}=\zeta$, where $\zeta $ is a
primitive $m^{\,\mathrm{th}}$ root of unity, $m\in \N $. Then we have
canonical isomorphisms
\begin{align*}
\deR (\gd)&\cong \C[\dett ^m,\dett ^{-m}]\ot \deR (\gdi) &&\text{for $\glqn$,}
\\
\deR (\gd)&\cong \deR (\gdi)\oplus \dett \deR (\gdi) &&\text{for $\oqn[2n+1]$,
$\Gamm=\Gamm_{\tau,1}$,}
\\
\deR (\gd)&\cong \deR (\gdi)\oplus  \deR (\dett \gdi)
&&\text{for $\oqn[2n+1]$,}
\\
\deR (\gdl)&\cong \deR (\gdr)\cong \deR (\gdi) &&\text{in  all cases.}
\end{align*}
{\em (c)} Let $G_q=\oqn[2n]$, $\A=\OO(G_q)$, and $\Gamm$ as in
{\em (a)}. Then  we have
\begin{align*}
\deR (\gd)\cong &\deR (\gdi)\oplus \dett  \deR (\gdi)\oplus\\
&\oplus \deR (\glm[\wedge]{\bo 1,\dett})\oplus \deR 
( \glm[\wedge]{\dett, \bo 1}),
\\
\deR (\gdl)\cong &\deR (\gdi)\oplus \deR (\glm[\wedge]{\bo 1,\dett}),
\\
\deR (\gdr)\cong &\deR (\gdi)\oplus \deR (\glm[\wedge]{\dett,\bo 1}).
\end{align*}
\end{theorem}
\begin{remark}
For $\glqn$ the case $z=1$ is of special interest since the commutation
relations between the differentials $\dd u^i_j$ and the matrix elements
$u^m_n$ appear to be {\em linear}, i.\,e.~
\begin{align}
\dd \bo u_1\cdot\bo u_2&=\Rda[\tau]{}\bo u_1\cdot\dd\bo u_2\Rda[\tau]{}.
\end{align}
These calculi were extensively studied in \cite{a-Sudbery2},
\cite{a-Malt1} and \cite{a-Tsygan1}. In the proof of Lemma\,\ref{l-nonzero}
given below we will show that $z=1$ is regular.
\end{remark}

\begin{theorem}\label{t-hodge}
Suppose that $q$ is a transcendental complex number.
Let $G_q$ be one of the quantum groups $\glqn$ or $\slqn$, and $\A=\OO(G_q)$
its coordinate Hopf algebra. Let $\Gamm$ be one of
the $N^2$-dimensional bicovariant first order differential calculi $\gtz$
over $\A $, where $z$ is regular in case of
$\glqn$. Then we have the Hodge decompositions
\begin{align}\label{e-hodge}
\begin{split}
\gd[k]&\cong \dd\,\gd[k-1]\oplus \kod[+]\,\gd[k+1]\oplus \deR ^k(\gd),
\\
\gd[k]&\cong \dd\,\gd[k-1]\oplus \kod[-]\,\gd[k+1]\oplus \deR ^k(\gd)
\end{split}
\end{align}
for $k\in\NO$. 
Moreover, the cohomology ring of $\gd $ is isomorphic to the algebra of
coinvariant forms and to the vector space of harmonic forms:
\begin{align} \label{e-cohom}
\deR ^k(\gd)\cong \deR ^k(\gdi)\cong \gdi[k]
\cong \har ^+ (\gd[k])\cong \har ^- (\gd[k]).
\end{align}
\end{theorem}

\begin{remarks}
(i) In \cite[Theorem\,3.2]{a-Schueler1} it was shown that $\gdi$ is a graded
commutative algebra and its Poincar\' e series has the form
$(1+t)(1+t^3)\cdots (1+t^{2N-1})$.
By the above theorem, $\dim \deR ^{N^2}(\gd)=\dim \gdi[N^2]=1$. This means
that there exists a linear functional $f$ on $\A$ with the following property.
For all $a\in \A $ there exists $\rho\in \gd[(N^2-1)]$ such that
$$
a\nu=\dd  \rho+f(a)\nu,
$$
where $\nu\in\gdi[N^2]$ is the unique up to scalars coinvariant form of degree
$N^2$ (volume form).
{}From the fact that $\dd $ is bicovariant, one derives easily that
$a_{(1)}f(a_{(2)})=f(a_{(1)})a_{(2)}=f(a)1$ and $f(1)=1$. Therefore,
$f$ is the Haar functional $h$ of the cosemisimple Hopf algebra $\A$. 

(ii) If $\A$ belongs to the B-, C-, or D-series, then a coinvariant form is
not closed in
general. However, there is a weaker form of the decomposition \rf[e-hodge].
Let $\Lam^k:=\Lapp[](\gd[k])=\Lapm[](\gd[k])$. Then
one can prove that
\begin{align}
\dd\, \Lam^{k-1}\oplus \kod[+]\,\Lam^{k+1}\cong\Lam^k&\cong \dd\,
\Lam^{k-1}\oplus \kod[-]\,\Lam^{k+1}.
\end{align}

Using a computer algebra program we calculated the first terms of the
Poincar\' e series $P(\gdi,t)=1+t+5t^3+15t^4+\cdots$ and
$P(\deR (\gdi),t)=1+t+2t^3+2t^4+\cdots$.
\end{remarks}

\section{Duality of Hopf bimodules}\label{duality}

We will show that the notion of  a $\sig$-metric naturally emerges by
considering the left-dual and right-dual Hopf bimodules of a given
Hopf bimodule. 
Our main result states that
the left-dual Hopf bimodule is isomorphic to the right-dual
Hopf bimodule.
This makes the  notion  of a bicovariant \mbox{$\sig$-metric} more
transparent. However, both notions are not identical since
the dual pairings $\ldual$ and $\rdual$ are not completely
\mbox{$\sig$-symmetric} while  the metric $g$ is.
 
\begin{definition}\label{d-dual}
Suppose that $\Gamm$ is a Hopf bimodule. A Hopf bimodule $\veeg$ is called
the {\em left-dual} to the Hopf bimodule $\Gamm$ if there exists a
homomorphism $\ldual\colon\veeg\ota\Gamm\to\A$ of Hopf bimodules such that
the pairing $\ldual$  is non-degenerate. Similarly, a Hopf bimodule $\gvee$
is called the {\em right-dual} to the Hopf bimodule $\Gamm$ if there exists
a non-degenerate homomorphism $\rdual\colon\Gamm\ota\gvee\to\A$ of Hopf
bimodules.
\end{definition}

\begin{remarks} (i) Note that the left-dual $\veeg$ and the right dual $\gvee$
to the Hopf bimodule $\Gamm $ always exist. Moreover, they are unique
up to isomorphisms. Indeed,
the projection $P_\ll\,\rho:=S(\rho_{(1)})\rho_{(0)}$ onto the
left-coinvariant subspace commutes with $\ldual$. Hence $\veeg_\ll$ and
$\gl$ are dually paired vector spaces.
Suppose that $\Gamm=(\bo v,\bo f)$.
It was shown in \cite[Section~II.\,4]{a-BCDRV} that
$\veeg=(\bo v^\fettc,\f)$. Similarly one proves that
$\gvee=({^\fettc\bo v},\bo f^\fettc)$. 

(ii) Since $\ldual$ and  $\rdual$ are homomorphisms of bicomodules, they are
bicovariant, i.\,e.~
$(\id\ot\ldual)\dl=\Delta\ldual=(\ldual\ot\id)\dr$ and
$(\id\ot\rdual)\dl=\Delta\rdual=(\rdual\ot\id)\dr$
on $\veeg\ota\Gamm$ and $\Gamm\ota\gvee$, respectively. 

(iii) It was shown in  \cite[Section~II.\,4]{a-BCDRV} that the pairing is
compatible with the braiding $\sig$. More precisely, let $\Lam$ be a Hopf 
bimodule. Then we have
$\ldual[23]\sig_{12}\sig_{23}=\ldual[12]$ on $\veeg\ota\Gamm\ota\Lam$
and $\ldual[12]\sig_{23}\sig_{12}=\ldual[23]$ on
$\Lam\ota\veeg\ota\Gamm$. Similarly, $\rdual[23]\sig_{12}\sig_{23}=\rdual[12]$
on $\Gamm\ota\gvee\ota\Lam$ and $\rdual[12]\sig_{23}\sig_{12}=\rdual[23]$ on
$\Lam\ota\Gamm\ota\gvee$.
\end{remarks}
\begin{proposition}\label{p-equiv}
Let $\A$ be a Hopf algebra with invertible antipode and let
$\Gamm=(\bo v,\bo f)$ be a Hopf bimodule with basis $\{\om_i\}$ of $\gl$. 
Let   $\{\theta_i\}$ and $\{\eta_i\}$ denote the left-coinvariant bases of
$\veeg$ and $\gvee$, dual to $\{\om_i\}$, respectively.
\\
Then the linear mapping $T\colon\veeg\to \gvee$ defined by
$T(a\,\theta_i)=a\,S^2(\eta_i)=a\,S(v^i_j)\eta_k\,v^j_k$
is an isomorphism of Hopf bimodules, where  $S$ denotes the antipode in the
graded super Hopf algebra
$(\gvee)^\ot$. 
Moreover, we have the following $\sig$-symmetry of the above pairing:
$$
\rdual(\id\ot T)=\ldual\sig\qquad\text{ on }\quad \Gamm\ota\veeg.
$$
\end{proposition}

\begin{proof}
(a) We first show that $T$ is a right comodule mapping, i.\,e.\ 
$\dr(\tilde{\theta}_i)=\tilde{\theta }_j\ot (\bo v^\fettc)^j_i$ for
$\tilde{\theta}_i:=T(\theta_i)=T^j_i\,\eta_j$.
Recall that the coproduct on $\gvee$ is given by
$\Delta(\eta_i)=\dl(\eta_i)+\dr(\eta_i)=1\ot\eta_i+\eta_j\ot
S^{-1}(v^i_j)$, see \cite[Proposition~13.7]{b-KS}. Since $\ve(\eta_i)=0$ one has $S(\eta_i)=-\eta_jv^i_j$ and
consequently, 
\begin{equation}\label{e-theta}
\tilde{\theta_i}=S^2(\eta_i)=-Sv^i_j\,S(\eta_j)=Sv^i_j\,\eta_k\,v^j_k
=\eta_k\rac v^i_k. 
\end{equation}
Hence
$\dr \,\tilde{\theta_i}=S v^x_j\,\eta_z\,v^j_y\ot S v^i_x\,S^{-1} v^k_z \,v^y_k=\tilde{\theta_x}\ot
Sv^i_x$ which proves (a).

(b) We show that $T$ is a right module map. Since $S^2$ is an algebra map of
$(\gvee)^\otimes$ we have for $a\in\A$
\begin{align*}
\tilde{\theta_i}\,S^2(a)&=S^2(\eta_i\, a)=S^2\bigl((S(f^j_i)\ast a)\,\eta_j\bigr)
\\
&=S^2\bigl(a_{(1)}S f^j_i(a_{(2)})\bigr)\tilde{\theta_j}=(S^2
a)_{(1)}S^{-1}f^j_i\bigl( (S^2 a)_{(2)}\bigr)\tilde{\theta_j}
\\
&=(S^{-1} f^j_i\ast S^2 a)\tilde{\theta_j}=T\bigl( (S^{-1} f^j_i\ast S^2 a)\,\theta_j\bigr)
\\
&=T(\theta_i\, S^2(a)).
\end{align*} 
Since $S^2\colon\A\to\A$ is surjective, $T$ is a right module map. By
\rf[e-theta], $\tilde{\theta_i}$ is a left-coinvariant 1-form. Hence $T$ 
is a left comodule map. Therefore $T$ is an isomorphism of Hopf 
bimodules. In particular, $\overline{T}:=(T^l_i)
\in\Mor(\bo v^\fettc,{^\fettc\bo v})$ and
$\overline{T}^\fettt\in\Mor(\bo f^\fettc, \f)$. 

(c) We prove the last assertion. Since $\bo f^\fettc$ defines the right action
on $\eta_i$, we obtain from \rf[e-theta] that
$\tilde{\theta_i}=Sv^i_j(S f^l_k\ast v^j_k)\eta_l=S f^l_k(v^i_k)\eta_l$,
i.\,e.\ $T^l_i=Sf^l_k(v^i_k)$. 
Since $\sig$, $g_\ll$, $g_\rr$, and $T$ are bimodule maps it suffices to prove
the statement for $\om_i\in\gl$ and $\theta_j\in\veeg_\ll$. By \rf[e-sigli] 
\begin{align*}
g_\ll\sig(\om_i\ota\,\theta_j)&=\ldual(\theta_k\ota (\om_i\rac S v^j_k))
=\ldual(\theta_k\ota \,f^i_n(S v^j_k)\,\om_n)
\\
&=f^i_n(Sv^j_n)=T^i_j=\rdual( \om_i\ota  T(\theta_j)).
\end{align*}
\end{proof}
\begin{remarks} (i)
Unfortunately the equation ${\ldual\, (T^{-1}\ot \id)}=\rdual\,\sig$
on $\gvee\ota\Gamm$ is not fulfilled in general. If this symmetry holds,
then the matrices $\overline{T}$ and $\tilde{T}=(\tilde{T}^a_b)$,
$\tilde{T}^a_b=S f^k_b(v^k_a)$, have to be inverse to each other. This
is not the case for the fundamental Hopf bimodules $(\bo u,\bo
\ell^{\pm \fettc})$, but for the differential Hopf bimodules $\gtz$
it is. A sufficient condition for the second $\sig$-symmetry is that
$\Gamm=(\bo v,\bo f)$ is an irreducible Hopf bimodule and that both
$\gi$ and $\Gamm^\vee_\ii$ are nontrivial.

(ii) It is easy to show that ${}^\vee\!\gpz\cong\gmz$.
Moreover,
the $\sig$-metric $g\colon\gpz\ota\gmz\to\A$, see \rf[e-metric], can be
obtained from $\rdual$ by the above identification of $\gmz$ with the
right-dual Hopf bimodule $\Gamm _{+,z}^\vee $ of $\gpz $:
\begin{equation}\nn
g(\om^+_{ij},\om^-_{kl})=g(\om^+_{ij},T^{-1}(\om_{kl}^{\vee})):=
g_\rr(\om^+_{ij},(T^{-1})^{mn}_{kl}\om^\vee_{mn})=(T^{-1})^{ji}_{kl},
\end{equation}
where $\{\om^\vee_{kl}\}$ is the left-coinvariant basis of $\Gamm ^\vee _{+,z}$
such that $g_\rr (\om ^+_{ij},\om ^\vee _{mn})=\delta _{jm}\delta _{in}$
and $T^{mn}_{kl}=(\ell^{-\fettc}\ot\ell^+)^{mn}_{xy}
\bigl(S(u^k_x\,(u^\fettc{})^l_y)\bigr)$ by step (c) of the above proof.
\end{remarks}

\section{Properties of the contraction}

We summarise some useful properties of the contraction, see
\cite[Lemmata~4.3, 4.4, 6.2]{a-Heck1}.
For $\xi _i\in \Gaw[k_i]{\tau _i}$, $i=0,1,2$, $\tau _1=\tau _2=-\tau _0$,
$k_1+k_2\leq k_0$, the contractions satisfy the following relations:
\begin{equation}\label{e-ctrlr}
\begin{split}
\ctrpm{\xi _1}{\ctrpm{\xi _2}{\xi _0}} &=
\ctrpm{\xi _1 \wedge \xi _2}{\xi _0},\,
\ctrpm{\ctrpm{\xi _0}{\xi _1}}{\xi _2}=
\ctrpm{\xi _0}{\xi _1\wedge \xi _2},
\\
\ctrpm{\xi _1}{\ctrpm{\xi _0}{\xi _2}} &=
\ctrpm{\ctrpm{\xi _1}{\xi _0}}{\xi _2}. 
\end{split}
\end{equation}
For  $a\in \A$, $\rho \in (\Gamm_\tau)_\ll$, and 
$\zeta\in (\Gamm_{-\tau})_\ll$ we have $h\ctrpm{a\rho }{\zeta}=h\ctrpm{\rho
  a}{\zeta} =h(a)\ctrpm{\rho }{\zeta}$. Since  $\ctrpm{\cdot}{\cdot}$ is
$\sig$-symmetric, for $\xi\in\gt$ we particularly get
\begin{gather}\label{e-haar} 
h\ctrpm{\xi}{\om_0^{-\tau}}=h\ctrpm{\om_0^{-\tau}}{\xi}.
\end{gather}
For  $\xi \in \Gaw[k]{\tau }$, $\xi '\in \Gaw[k]{-\tau }$,
$\rho _1\in \Ga{\tau}$, $\rho _2\in \Ga{-\tau }$, $k\geq 1$,  we have
\begin{align}
\label{eq-ctrrekr} 
\begin{split}
\ctrpm{\xi \wedge \rho _1}{\rho _2}&=\xi \ctrpm{\rho _1}{\rho _2}-
\ctrpm{\xi }{\rho ^\mp _{(1)}}\wedge \rho ^\mp _{(2)},\\ 
\ctrpm{\rho _1}{\rho _2\wedge \xi '}&=\ctrpm{\rho _1}{\rho _2}\xi '
-\rho ^\mp _{(1)}\wedge \ctrpm{\rho ^\mp _{(2)}}{\xi '},
\end{split}
\end{align}
 where
$\sigma ^\mp (\rho _1\ota \rho _2)=\rho ^\mp _{(1)}\ota \rho ^\mp _{(2)}
\in \Ga{-\tau }\ota \Ga{\tau }$.

Now let us prove an identity for the braiding morphism $\sig$.
\begin{lemma}\label{l-sig}
Let $\Gamm$ and $\Lam$ be Hopf bimodules over $\A$. Then  we have
\begin{align}\label{e-sig}
(\sig_{\Gamm,\Lam})_k\cdots(\sig_{\Gamm,\Lam})_1&=\sig_{\Gamm,\Lamb[\ot k]}
\end{align}
and this map is a homomorphism of the Hopf bimodule
$\Gamm\ota\Lam\ota\cdots\ota\Lam$ to $\Lam\ota\cdots\ota\Lam\ota\Gamm$. 
Moreover, replacing $\Lamb[\ot k]$ by its quotient $\Lamb[\land k]$,
equation \rf[e-sig] remains valid. Similarly,
$\sig^-_k\cdots\sig^-_1\colon\Gamm\ota\Lamb[\wedge k]\to\Lamb[\wedge
k]\ota\Gamm$ is well-defined and coincides with $\sig^-_{\Gamm,\Lamb[\wedge k]}$.
\end{lemma}
\begin{proof}
(a) The braiding $\sig$ is compatible with the tensor product of Hopf 
bimodules in the sense that the identity
${(\id_Y\ot\sig_{X,Z}){\kri}(\sig_{X,Y}\ot\id_Z)}={\sig_{X,Y\ot Z}}$ is fulfilled for all
Hopf bimodules $X$, $Y$, and $Z$, see \cite[Theorem~5.2]{a-Yetter}. Iterating
this yields
$$
(\sig_{\Gamm,\Lam})_k(\sig_{\Gamm,\Lam})_{k-1}\cdots(\sig_{\Gamm,\Lam})_{1}=\sig_{\Gamm,\Lamb[\ot
  k]}.
$$
(b) In what follows we skip the space indices $\Gamm$ and $\Lam$
to simplify the notations. Since
$\sig_{i+1}\sig_{i}\sig_{i+1}=\sig_i\sig_{i+1}\sig_{i}$
we obtain
$\sig_k\cdots\sig_1\,\sig_{i+1}=\sig_{i}\,\sig_k\cdots\sig_1$,
$i=1,\dots,k-1$. It
follows that
$\sig_k\cdots\sig_1\,(A_k)_{2\cdots
  k+1}=A_k\,\sig_k\cdots\sig_1$.
Hence $\sig_k\cdots\sig _1$ maps $\Gamm \ota \ker A_k$ to $\ker A_k\ota \Gamm $
and therefore it defines a mapping $\Gamm \ota \Lamb[\land k]\to
\Lamb[\land k]\ota \Gamm $. Moreover, by \rf[e-sig] it coincides with
$\sig_{\Gamm,\Lamb[\land k]}$.
Since $\sig^-$ defines a braiding as well and $\ker A_k=\ker A^-_k$,
the proof for $\sig^-$ is analogous.
\end{proof}

Now we  add some new relations which not yet appeared in \cite{a-Heck1}.
\begin{lemma}\label{l-ctrm}
For $\xi\in\Gaw[k]{\tau}$, $\rho_1\in\Ga{\tau}$, $\rho_2\in\Ga{-\tau}$, $k\ge
1$, $\tau\in\{+,-\}$, we
have
\begin{align}\label{e-ctr1}
\begin{split}
\ctrpm{\rho_1\wedge\xi}{\rho_2}&=\rho_1\wedge\ctrpm{\xi}{\rho_2}+(-1)^k\gtil(\sig^\pm(\rho_1\ota\xi),\rho_2),
\\
\ctrpm{\rho_2}{\xi\wedge\rho_1}&=\ctrpm{\rho_2}{\xi}\wedge\rho_1+(-1)^k\gtil(\rho_2,\sig^\pm(\xi\ota\rho_1)),
\end{split}
\end{align}
where $\sig^\pm=\sig^\pm_{\gt,\Gaw[k]{\tau}}$ in the first equation and
$\sig^\pm=\sig^\pm_{\Gaw[k]{\tau},\gt}$ in the second equation.
\end{lemma}
\begin{proof}
We carry out the proof of the first equation. The proof of the  second one
is analogous. By
definition \rf[eq-ctr],
\begin{align*}
\ctrpm{\rho_1\wedge\xi}{\rho_2}&=\gtil\bigl((1-\sig^\pm_{k}+\sig^\pm_{k}
\sig^\pm_{k-1}-\cdots+(-1)^{k}\sig^\pm_{k}\cdots\sig^\pm_{1})(\rho_1\ota\xi),
\rho_2\bigr).
\end{align*}
Note that the  endomorphisms $\sig^\pm_{k}$,
$\sig^\pm_{k}\sig^\pm_{k-1}$, \dots,
$\sig^\pm_{k}\cdots\sig^\pm_{2}$ do {\em not} act in the first component. So
we can separate the last summand.
Applying  Lemma \ref{l-sig} we continue
\begin{align*}
&=\gtil\bigl(\rho_1\ota(1-\sig^\pm_{k-1}+\cdots+(-1)^{k-1}\sig^\pm_{k-1}\cdots\sig^\pm_{1})\xi,\rho_2\bigr)
+
\\
&\quad +(-1)^k\gtil(\sig^\pm_{k}\cdots\sig^\pm_{1}(\rho_1\ota\xi),\rho_2)
\\
&=\rho_1\ota\ctrpm{\xi}{\rho_2}+(-1)^k\gtil(\sig^\pm(\rho_1\ota\xi),\rho_2).
\end{align*}
This finishes the proof.
\end{proof}
\section{Quantum Laplace-Beltrami Operators}\label{laplace}

Let us derive some important properties of the quantum Laplace-Beltrami
operators defined by \rf[e-laplace].
\begin{proposition}\label{p-lapsigma} For $\rho\in\Gaw[k]{\tau}$ we have
\begin{align}\label{e-lapsigma}
\Lappm\rho=(-1)^{k}\bigl(-2\sss\rho+\gtil\bigl(\sig^\pm(\om^\tau_0\ota\rho),\om^{-\tau}_0\bigr)+\gtil\bigl(\om^{-\tau}_0,\sig^\pm(\rho\ota\om^\tau_0)\bigr)\bigr).
\end{align}
\end{proposition}
\begin{proof}
Using the definition of $\dd_\tau$ and $\kodpm_\tau$,
the identity $\sig(\om^{-\tau}_0\ota\om^\tau_0)=\om^\tau_0\ota\om^{-\tau}_0$
and equations \rf[eq-ctrrekr] and \rf[e-ctrom] we get
(the not underlined terms remain unchanged)
\begin{align*}
\Lappm\rho&=-\dd_\tau \bigl(\ctrpm{\rho}{\om^{-\tau}_0}+(-1)^k\ctrpm{\om^{-\tau}_0}{\rho}\bigr)+\kodpm_\tau\bigl(\om^\tau_0\wedge\rho+(-1)^{k-1}\rho\wedge\om^\tau_0\bigr)
\\
&=-\om^\tau_0\wedge\ctrpm{\rho}{\om^{-\tau}_0}
+(-1)^{k-1}\ctrpm{\rho}{\om^{-\tau}_0}\wedge\om^\tau_0+(-1)^{k+1}\bigl(\om^\tau_0\wedge\ctrpm{\om^{-\tau}_0}{\rho}
\\
&\quad+(-1)^k\ctrpm{\om^{-\tau}_0}{\rho}\wedge\om^\tau_0\bigr)+
\ctrpm{\om^\tau_0\land\rho}{\om^{-\tau}_0}+\ul{(-1)^{k+1}\ctrpm{\om^{-\tau}_0}{\om^\tau_0\land\rho}}
\\
&\quad+(-1)^{k-1}(\ul{\ctrpm{\rho\wedge\om^\tau_0}{\om^{-\tau}_0}}+(-1)^{k+1}\ctrpm{\om^{-\tau}_0}{\rho\wedge\om^\tau_0})
\\
&=-\om^\tau_0\wedge\ctrpm{\rho}{\om^{-\tau}_0}
+\ul{(-1)^{k-1}\ctrpm{\rho}{\om^{-\tau}_0}\wedge\om^\tau_0}+\ul{(-1)^{k+1}\om^\tau_0\wedge\ctrpm{\om^{-\tau}_0}{\rho}}
\\
&\quad-\ctrpm{\om^{-\tau}_0}{\rho}\wedge\om^\tau_0+
\ctrpm{\om^\tau_0\land\rho}{\om^{-\tau}_0}+(-1)^{k+1}(\sss\rho\ul{-\om^\tau_0\wedge\ctrpm{\om^{-\tau}_0}{\rho}})
\\
&\quad +(-1)^{k-1}(\sss\rho\ul{-\ctrpm{\rho}{\om^{-\tau}_0}\wedge\om^\tau_0})+\ctrpm{\om^{-\tau}_0}{\rho\wedge\om^\tau_0}
\\
&=2\sss(-1)^{k+1}\rho-\om^\tau_0\wedge\ctrpm{\rho}{\om^{-\tau}}-\ctrpm{\om^{-\tau}_0}{\rho}\wedge\om^\tau_0+\ul{\ctrpm{\om^\tau_0\land\rho}{\om^{-\tau}_0}}+
\\
&\quad +\ul{\ctrpm{\om^{-\tau}_0}{\rho\wedge\om^\tau_0}}
\\
&=2\sss(-1)^{k+1}\rho+(-1)^k\bigl(\gtil(\sig^\pm(\om^\tau_0\ota\rho),\om^{-\tau}_0)+\gtil(\om^{-\tau}_0,\sig^\pm(\rho\ota\om^\tau_0))\bigr).
\end{align*}
The last equation follows by \rf[e-ctr1].
\end{proof}
\begin{lemma}\label{l-00}
For $a\in\A$ and $\tau\in\{+,-\}$ we have
\begin{equation}\label{e-00}
g(\om^\tau(a),\om^{-\tau}_0)=g(\om^{-\tau}_0,\om^\tau(a))=\sss X^\tau_0(a).
\end{equation}
\end{lemma}
\begin{proof}
Let $\om^\tau_{1,ij}:=\om^\tau_{ij}-(P_0^\tau)^{kl}_{ij}\om^\tau_{kl}$,
where $P^\tau_0$ is given by \rf[e-P0].
Further let $X^{\tau}_{1,ij}$  be the corresponding dual
basis elements in the quantum tangent space of $\gt$.
Since the $\sig$-metric is bicovariant,
the complex matrices
$(g(\om^\tau_{1,ij},\om^{-\tau}_0))$ and $(g(\om^{-\tau}_0,\om^\tau_{1,ij}))$
are elements of the vector space $\Mor(\bo 1,\bo v\ot \bo 1)$,
where $\bo v$ denotes the corepresentation corresponding to the right
coaction on $\langle \om ^\tau _{1,ij} \rangle $.
By Schur's lemma these matrices have to be zero 
(see also the proof of Proposition \ref{p-nondeg}\,(i)).
Using \rf[e-ctrom] and
$\om^\tau(a)=X^\tau_0(a)\om^\tau_0+\sum_{ij}X^\tau_{1,ij}(a)\om^\tau_{1,ij}$
the assertion follows.
\end{proof}
\begin{proposition}\label{p-ev}
For $\rho\in\Gaw[k]{\tau}$, $\tau\in\{+,-\}$, we have
\begin{align}
\begin{split}
\Lapp\rho&=
(-1)^k\sss\bigl(\rho_{(0)}X^\tau_0(\rho_{(1)})+X^{-\tau}_0(\rho_{(-1)})\rho_{(0)}\bigr)
\\
&=(-1)^k\sss\bigl(X_0^\tau\ast\rho+\rho\ast X_0^{-\tau}\bigr),
\end{split}
\label{e-lp}\\
\begin{split}
\Lapm\rho
&=(-1)^k\sss\bigl(\rho_{(0)}X^{-\tau}_0(\rho_{(1)})+X^\tau_0(\rho_{(-1)})\rho_{(0)}\bigr)
\\
&=(-1)^{k}\sss\bigl(X_0^{-\tau}\ast\rho+\rho\ast X_0^\tau\bigr).
\end{split}
 \label{e-lm}
\end{align}
In particular, $\Lapp a=\Lapm  a=\sss X_0\ast a=\sss a\ast X_0$ for $a\in\A$.
\end{proposition}
\begin{proof}
We prove \rf[e-lp]. Let
$\rho=\sum _i\rho_ia_i$  be a presentation of $\rho$ with
$\rho_i\in(\gdt[k])_\ll$, $a_i\in\A$. By \rf[e-sigli], since $g$ and $\sig$ are  $\A$-module homomorphisms and
since $\beta a=a_{(1)}\beta\rac a_{(2)}$, 
 we obtain
\begin{align*}
\gtil(\sig(\om^\tau_0\ota\rho_ia_i),\om^{-\tau}_0)&=\gtil(\sig(\om^\tau_0\ota\rho_i)a_i,\om^{-\tau}_0)
\\
&=\gtil(\rho_i{}_{(0)}\ota(\om^\tau_0\rac\rho_i{}_{(1)})a_i,\om^{-\tau}_0)
\\
&=\gtil(\rho_i{}_{(0)}a_i{}_{(1)}\ota(\om^\tau_0\rac(\rho_i{}_{(1)}a_i{}_{(2)})),\om^{-\tau}_0)
\\
&=\rho_i{}_{(0)}a_i{}_{(1)}g\bigl(\om^\tau(\rho_i{}_{(1)}a_i{}_{(2)})+\ve(\rho_i{}_{(1)}a_i{}_{(2)})\om^\tau_0,\om^{-\tau}_0\bigr)
\\
&=\rho_i{}_{(0)}a_i{}_{(1)}X^\tau_0(\rho_i{}_{(1)}a_i{}_{(2)})+\rho_ia_ig(\om^\tau_0,\om^{-\tau}_0)
\\
&=\sss(\rho_{(0)}X^\tau_0(\rho_{(1)})+\rho).
\end{align*}
In the fourth equation we used \rf[e-oma], in the fifth equation \rf[e-00] and in the last one \rf[e-ctrom].
Now let $\rho=\sum _ia_i\rho_i$, $\rho_i\in(\Gaw[k]{\tau})_\ll$,
$a_i\in\A$. By \rf[e-00] and \rf[e-ctrom] we have
\begin{align*}
\gtil(\om^{-\tau}_0,\sig(a_i\rho_i\ota\om^\tau_0))&=\gtil(\om^{-\tau}_0a_i,\sig(\rho_i\ota\om^\tau_0))
\\
&=\gtil(a_i{}_{(1)}\om^{-\tau}_0\rac
a_i{}_{(2)},\om^\tau_0\ota\rho_i)
\\
&=a_i{}_{(1)}g\bigl(\om^{-\tau}(a_i{}_{(2)})+\ve(a_i{}_{(2)})\om^{-\tau}_0,\om^\tau_0\bigr)\rho_i
\\
&=\sss\bigl(a_i{}_{(1)}X^{-\tau}_0(a_i{}_{(2)})\rho_i+a_i\rho_i\bigr)
\\
&=\sss(X^{-\tau}_0(\rho_{(-1)})\rho_{(0)} +\rho).
\end{align*}
In the last equation we used
$b_{(1)}X^{-\tau}_0(b_{(2)})=X^{-\tau}_0(b_{(1)})b_{(2)}$ for all
$b\in\A$ which follows from the centrality of $X^{-\tau}_0$. 
Inserting both parts  into \rf[e-lapsigma] we obtain \rf[e-lp].

Let us prove \rf[e-lm]. Similarly to the preceding equation one
shows that $\gtil(\sig^-(\om^\tau_0\ota\rho),\om^{-\tau}_0)
=\sss\bigl(\rho+X_0(\rho_{(-1)})\rho_{(0)}\bigr)$. 
Let $\rho=\sum _ia_i\rho_i$ with left-coinvariant elements $\rho_i$. Using
\rf[e-sigli] we get
\begin{align*}
\gtil(\om^{-\tau}_0,&\sig^-(a_i\rho_i\ota\om^\tau_0))=\gtil(\om^{-\tau}_0a_i,(\om^\tau_0\rac S^{-1}\rho_i{}_{(1)})\ota\rho_i{}_{(0)})
\\
&=g\bigl(\om^{-\tau}_0a_iS(S^{-1}\rho_i{}_{(2)}),\om^\tau_0S^{-1}\rho_i{}_{(1)}\bigr)\rho_i{}_{(0)}
\\
&=g\bigl(a_i{}_{(1)}\rho_i{}_{(2)}\om^{-\tau}\rac(a_i{}_{(2)}\rho_i{}_{(3)}),\om^\tau_0\bigr)S^{-1}\rho_i{}_{(1)}\rho_i{}_{(0)}
\\
&=\sss a_i{}_{(1)}\rho_i{}_{(2)}\bigl(\ve(a_i{}_{(2)}\rho_i{}_{(3)})+X^{-\tau}_0(a_i{}_{(2)}\rho_i{}_{(3)})\bigl)S^{-1}\rho_i{}_{(1)}\rho_i{}_{(0)}
\\
&=\sss(\rho+a_i{}_{(1)}\rho_i{}_{(0)}X^{-\tau}_0(a_i{}_{(2)}\rho_i{}_{(1)}))
\\
&=\sss(\rho+\rho_{(0)}X^{-\tau}_0(\rho_{(1)})).
\end{align*}
This gives \rf[e-lm].  Since $X_0^\tau$  is  central,
$\Lapp  a=\Lapm  a=\sss X_0\ast a=\sss a\ast X_0$
follows from \rf[e-lp] and \rf[e-lm].
\end{proof}
\subsection*{Cosemisimple Hopf algebras}
Let $\A$ be a cosemisimple Hopf algebra (cf. \cite[Subsection~11.2]{b-KS}) and
let  $\AD$ be the set of equivalence classes $\alpha$ 
of irreducible corepresentations $\bo u^\alpha$ of $\A$. Then $\A$ has  the  Peter-Weyl decomposition 
$\A=\bigoplus_{\alpha\in\AD}\CC(\bo u^\alpha)$. Let $P^\alpha\colon\A\to\A$
denote the projection of $\A$ onto the simple coalgebra $\CC(\bo u^\alpha)$. In
particular, if $a=\sum_{\lam,i,j}c^\lam_{ij}u^\lam_{ij}$, $c^\lam_{ij}\in\C$,
and $\{u^\lam_{ij}|\, i,\,j=1,\dots,d_\lam\}$ is a linear basis of
$\CC(\bo u^\lam)$, then $P^\alpha(a)=\sum_{ij}c^\alpha_{ij}u^\alpha_{ij}$. 
Define the linear functionals $h^\alpha$, $\alpha\in\AD$,  on $\A$ by
$h^\alpha=\ve{\kri} P^\alpha$. Obviously, we have
$\sum_{\alpha\in\AD}h^\alpha=\ve$
and $(P^\alpha\ot\id)\Delta=(\id\ot P^\alpha)\Delta=(P^\alpha\ot
P^\alpha)\Delta =\Delta{\kri} P^\alpha$. 
It is easily seen  that $h^\alpha\ast a=a\ast h^\alpha=P^\alpha(a)$. Note that
$h^0$ corresponding to  $\bo u^0\equiv\bo 1$ is the Haar functional on $\A$.

Since $X^\tau_0$ is central, for $\lam\in\AD$ there exist complex numbers $E^\tau_\lam$ such
that $X^\tau_0\ast$ acts as a scalar on $\CC(\bo u^\lam)$:
\begin{equation}\label{e-x}
\begin{split}
\sss X_0^\tau\ast h^\lam\ast a&=E^\tau_\lam h^\lam \ast a,
\\
\sss X_0^\tau(h^\lam\ast a)&=E^\tau_\lam h^\lam(a)
\end{split}
\end{equation} 
for $a\in\A$. Let $\bo v,\bo w$ be corepresentations of $\A$ and $\Gamm=\gt$.  Define the following subspaces of $\gd$:
\begin{align*}
\glm{\bo v,\bo w}&:=\{\rho\in\gd[k]\,|\, \dlr\in\CC(\bo v)\ot
\Gaw[k]{}\ot\CC(\bo w)\},
\\
\glm[\wedge]{\bo v,\bo w}&:=\bigoplus_{k\ge0}\glm{\bo v,\bo w}.
\end{align*}
We briefly write $\glm{\lam,\mu}$  instead of
$\glm{\bo u^\lam,\bo u^\mu}$. The main step in our proof is the following
spectral decomposition of the quantum Laplace-Beltrami operators.
\begin{proposition}\label{p-ev1}
Let $\Gamm=\gt$, $\tau\in\{+,-\}$. For $\lam,\,\mu\in\AD$ define the mapping
$h^{\lam\mu}\colon\gd\to\glm[\wedge]{\lam,\mu}$ by
$h^{\lam\mu}(\rho)=h^\lam(\rho_{(-1)})\rho_{(0)}h^\mu(\rho_{(1)})$.  
For  $\rho\in\gd$ and  $\rho^{\lam\mu}=h^{\lam\mu}(\rho)$
we then have
\begin{align} 
\rho&=\sum_{\lam,\mu\in\AD}\rho^{\lam\mu},\quad\gd=\bigoplus_{k\ge0}\bigoplus_{\lam,\mu\in\AD}\glm{\lam,\mu},\label{e-dec}
\\
\begin{split}  \label{e-lrho}
  \Lapp\rho^{\lam\mu}&=(-1)^k(E^{-\tau}_\lam+E^{\tau}_\mu)\rho^{\lam\mu},
   \\
   \Lapm\rho^{\lam\mu}&=(-1)^k(E^\tau_\lam +E^{-\tau}_\mu)\rho^{\lam\mu}.
\end{split}
\end{align}
For brevity we write $E_{\lam\mu}:=E^-_\lam+E^+_\mu$. 
\end{proposition}
\begin{proof} (a) An easy computation shows that indeed
   $\rho^{\lam\mu}\in\glm[\wedge]{\lam,\mu}$. 
Since $\sum_\lam h^\lam=\ve$ and
$\rho=\ve(\rho_{(-1)})\rho_{(0)}\ve(\rho_{(1)})$, the first part of
\rf[e-dec] follows. Let us verify the second part of \rf[e-dec].
The first sum is direct by the grading. The second sum is
direct, since matrix elements of inequivalent
irreducible corepresentations are linearly independent.

(b) Since $\dr(\rho^{\lam\mu})=
{h^\lam(\rho_{(-1)})\rho_{(0)}}\ot {h^\mu{\ast}\rho_{(1)}}$ and 
$\dl(\rho^{\lam\mu})=
{\rho_{(-1)}{\ast} h^\lam}\ot{ \rho_{(0)}h^\mu(\rho_{(1)})}$,
by \rf[e-lp] and \rf[e-x] we obtain the equation
\begin{align*}
\Lapp\rho^{\lam\mu}=&
(-1)^k\sss\bigl(h^\lam(\rho_{(-1)})\rho_{(0)}X^\tau_0(h^\mu\ast\rho_{(1)})+\\
&+X^{-\tau}_0(h^\lam\ast\rho_{(-1)})\rho_{(0)}h^\mu(\rho_{(1)})\bigr)
=(-1)^k(E^\tau_\mu+E^{-\tau}_\lam)\rho^{\lam\mu}.
\end{align*}
The proof of the second part of \rf[e-lrho] is analogous.
\end{proof}

In the remainder of this section $\A$ denotes the coordinate Hopf algebra of
the quantum group $G_q$ as in Theorem\,\ref{t-co}. For $\lam\in P_{++}$
as usual $\lam'_i$ denotes the length of the $i^{\,\mathrm{th}}$-column of
$\lam$. We define $P(\A )$ to be the set
$P_+$ for $\glqn$, $P_{++}^0:=\{\lam\in P_{++}\,|\,\lam_N=0\}$ for $\slqn$,
$\{\lam\in P_{++}\,|\,\lam'_1+\lam'_2\le N\}$ for $\oqn$,
$\{\lam\in P_{++}\,|\,\lam'_1\le n\}$ for $\spqn[2n]$,
and $\{\lam\in P_{++}\,|\,\lam'_1\le \frac{N}{2}\}$ for $\soqn$, respectively.
By \cite[Theorem~11.22]{b-KS} irreducible corepresentations
$\bo v$ of $\A$ are in one-to-one correspondence with elements of $P(\A )$.
We identify $\AD $ and $P(\A )$.

\begin{lemma}\label{l-nonzero}
Suppose that $\lam ,\mu \in P(\A )$.
\\
{\em (i)} For  $\slqn$  we have $E_{\lam\mu}=0$ if and
only if  $\lam=\mu=(0)$.
\\
{\em (ii)} 
For  $\glqn$ we have $E_{\lam\mu}=F_{\lam\mu}$.  The
parameter value $z=1$ is regular. If $z^Nq^{-2}=\zeta $ for a primitive
$m^{\,\mathrm{th}}$ root of unity $\zeta $, $m\in \N $,
then we have $E_{\lam\mu}=0$ if
and only if $\lam=(n^N)$ and $\mu=(k^N)$ for some $n,k\in m\Z$. 
\\
{\em (iii)} In the cases  $\spqn$ and $\soqn$  we have $E_{\lam\mu}=0$ if and
only if  $\lam=\mu=(0)$. In the case $\oqn$ we have $E_{\lam\mu}=0$ if and
only if $\lam,\mu\in\{(0),(1^N)\}$. 
\end{lemma}
\begin{proof}  For $\tau\in\{+,-\}$ define the following rational functions of $t$ and $z$:
\begin{align}
e^\tau_\lam(t,z)&:=z^{-\tau m}\bigl(\qnum_t
+\tau(t-t^{-1})\sum_{x\in\lam}\sgn(x)t^{\tau(N+2c(x))}\bigr)-\qnum_t,\notag
\\
\label{e-elampm}
e_{\lam\mu}(t,z)&:=e^-_\lam(t,z)+e^+_\mu(t,z),
\end{align}
where $m=|\lam|$.
It follows from \cite[Proposition~7.1]{a-Heck1} that for the quantum groups
$\glqn$ and $\slqn$ and for $\lam\in P_{++}$ we have
$E_\lam^\tau=e^\tau_\lam(q,z)$, where $z^N=q^2$ in the $\slqn$ case and $z\ne0$
in the $\glqn$ case.
Note that we have to replace $z^2$ in
\cite{a-Heck1} by $z$ according to our definition of $\gtz$. Obviously,
$E_{\lam\mu}=e_{\lam\mu}(q,z)$ for $\lam,\mu\in P_{++}^0$. 
Later we will see that $E_\lam^\tau=e^\tau_\lam(q,z)$ for  $\lam\in P_+$
and not only for $\lam\in P_{++}$.

We prove (i).  Set $\wtil{E}^\tau_\lam=\lim_{t\to
  1}(t-t^{-1})^{-2}\,
e^\tau_\lam(t,\,t^{2{/}N})$.
In the remark to
Proposition~7.1 in \cite{a-Heck1} it was noted that
\begin{equation}\label{e-lim1}
\wtil{E}_\lam:=\wtil{E}^+_\lam+\wtil{E}^-_\lam=\sum_{i=1}^{N-1}\frac{(N-i)m_i}{N}\bigl(i(m_i+N)+2\sum_{j=1}^{i-1}j
m_j\bigr),
\end{equation}
where $m_i=\lam_i-\lam_{i+1}$, $i=1,\dots,N-1$. On the other hand,
computing the limit $\lim_{t\to 1}(t-t^{-1})^{-2}
e^\tau_\lam(t,t^{2{/}N})$ directly  from \rf[e-elampm] one gets
\begin{equation}\label{e-lim2}
\wtil{E}^+_\lam=\wtil{E}^-_\lam=\frac{1}{2N}(m N^2+2 c(\lam)N -m^2).
\end{equation}
Using the formulae $c(\lam)=n(\lam')-n(\lam)$, $n(\lam)=\sum_{i\ge 1}(i-1)\lam_i$, and
${n(\lam')}={\sum_{i\ge 1}\half\lam_i(\lam_i-1)}$ from  \cite[Chapter~1]{b-Macdonald}, we
obtain \rf[e-lim1] from \rf[e-lim2].
{}From \rf[e-lim1] it follows that $\wtil{E}_\lam^+\ge 0$ for
$\lam\in P^0_{++}$ and $\wtil{E}_\lam^+=0$ if and only if $\lam=(0)$.
Suppose that $E_{\lam\mu}=0$ for some $\lam,\mu\in P^0_{++}$.
Since $e_{\lam\mu}(t):=e_{\lam\mu}(t,t^{2{/}N})$ is
an algebraic function of $t$ and  $t=q$ is a transcendental root,
$e_{\lam\mu}(t)\equiv 0$. In particular
\begin{align}
\lim_{t\to 1}(t-t^{-1})^{-2}e_{\lam\mu}(t,t^{2{/}N})
&=\wtil{E}^-_\lam+\wtil{E}^+_\mu=0.
\end{align}
Hence $\lam=\mu=(0)$.

Let us prove (ii).  We will show that $E^\tau_\lam=e^\tau_\lam(q,z)$
for arbitrary $\lam\in P_+$. For this purpose we  prove  that $\dett^n \, a$
is an eigenvector for  $X_0^\tau\ast\,$ if $a$ is and we
compute the corresponding eigenvalue. Suppose that $\sss X_0^\tau\ast a=E^\tau a$ for a complex number
$E^\tau$. 
Since $\ell^+{}^i_j(\dett)=qx^{-N}\delta_{ij}$ and
$\ell^-{}^i_j(\dett)=q^{-1}y^N\delta_{ij}$ we have
\begin{equation}\label{e-dd}
\om^\tau(\dett)=(q^{2\tau}z^{-\tau N}-1)\om_0^\tau.
\end{equation}
Hence $\om^\tau_0\rac
\dett=q^{2\tau}z^{-\tau N}\om_0^\tau$. Acting from the right by $\deti$ gives
$\om_0^\tau\rac\deti=q^{-2\tau}z^{\tau N}\om_0^\tau$. For $n\in\Z$ we thus have
$\om_0^\tau\rac \dett^n=q^{2n\tau}z^{-nN\tau}\om_0^\tau$. 
Since  $\dett^n$ is grouplike, $\om^\tau(a)=\om^\tau_0\rac
a-\ve(a)\om_0^\tau$, and $\rho\,a=a_{(1)}\rho\rac  a_{(2)}$,
we obtain by \rf[e-00] the following formulae for $n\in\Z $:
\begin{align*}
\sss X_0^\tau\ast(\dett^n a)&=\sss \dett^n a_{(1)}X_0^\tau(\dett^n a_{(2)})
\\
&=\dett^n a_{(1)}g\bigl(\om^\tau(\dett^n a_{(2)}),\om_0^{-\tau}\bigr)
\\
&=\dett^n a_{(1)}g(\om_0^\tau\rac (\dett^n a_{(2)}),\om_0^{-\tau})-\dett^n a_{(1)}\ve(a_{(2)})g(\om_0^\tau,\om_0^{-\tau})
\\
&=\dett^n a_{(1)} g(q^{2n\tau}z^{-nN\tau}\om_0^\tau\rac a_{(2)},\om_0^{-\tau})
-\sss \dett^n a
\\
&=q^{2n\tau}z^{-nN\tau}\dett ^n
a_{(1)}\bigl(g(\om^\tau(a_{(2)}),\om_0^{-\tau})+\ve(a_{(2)})\sss \bigr)-\sss
\dett^n a
\\
&=\bigl(q^{2n\tau}z^{-nN\tau}(E^\tau+\sss)-\sss\bigr)\dett^n a.
\end{align*}
Since $\dett$ corresponds to the weight $(1^N)$ and $\sss=\qnum_q$ we have for
$\lam\in P_{+}$ 
\begin{align}\label{e-elap}
E_{\lam+(1^N)}^\tau+\qnum_q=q^{2\tau}z^{-N\tau}(E_\lam^\tau+\qnum_q).
\end{align}
Next we will show that for $\lam\in P_+$,
\begin{equation}\label{e-recursion} 
e^\tau_{\lam+(1^N)}(t,z)+\qnum_t=t^{2\tau}z^{-N\tau}(e^\tau_\lam(t,z)+\qnum_t).
\end{equation} 
For $\lam\in P_+$ set
$\hat{\lam}=\{{(i,j)}\in\lam+(1^N)\,|\,{(i,j-1)}\in\lam\}$.  
How to obtain  $\lam+(1^N)$ from $\lam$? In cases $\lam_1\ge\dots\ge\lam_k\ge0$
shift $\lam$ to the right by $1$, then add one box $(i,1)$ for each
$i=1,\dots,k$. In cases $0\ge\lam_{k+1}\ge\dots\ge\lam_N$ remove the box
$(i,0)$ and then shift the remainder by $1$ to the right. Since boxes with
$\lam_i<0$ have negative sign and $|\lam+(1^N)|=m+N$ we obtain by \rf[e-elampm]
\begin{align*}
e^\tau_{\lam+(1^N)}(t,z)=&z^{-\tau(m+N)}\bigl(\qnum_t
+\tau(t-t^{-1})\sum_{x\in\hat{\lam}}\sgn(x)\,
t^{\tau(N+2c(x))}+
\\
&+\tau(t-t^{-1})\sum_{i=1}^Nt^{\tau(N+2-2i)}\bigr)-\qnum_t
\\
=&z^{-\tau m-\tau N}\bigl(\qnum_t
+\tau(t-t^{-1})t^{2\tau}\sum_{x\in\lam}\sgn(x)\,
t^{\tau(N+2c(x))}+
\\
&+\tau(t-t^{-1})(t^{\tau N}+t^{\tau(N-2)}+\cdots+t^{\tau(-N+2)})\bigr)-\qnum_t
\\
=&z^{-\tau m-\tau N}\bigl((1-t^{2\tau})\qnum_t
+z^{\tau m}t^{2\tau}(e^\tau_\lam(t,z)+\qnum_t)+
\\
&+\tau(t-t^{-1})t^\tau\qnum_t\bigr)-\qnum_t
\\
=&t^{2\tau} z^{-\tau N}(e_\lam ^\tau(t,z)+\qnum_t)-\qnum_t.
\end{align*}
In the last line we used $\tau t^\tau(t-t^{-1})=t^{2\tau}-1$.
Since there exists $n\in\N$ such that $\lam+(n^N)\in P_{++}$ and
$E^\tau_{\lam+(n^N)}=e^\tau_{\lam+(n^N)}(q,z)$ , from \rf[e-elap] and
\rf[e-recursion] we  obtain $E^\tau_\lam=e^\tau_\lam(q,z)$.
Comparing \rf[e-elampm] and \rf[e-hlam] yields $E_{\lam\mu}=F_{\lam\mu}$.

We will show that $z=1$ is regular.
Suppose that $F_{\lam\mu}=0$ for some $\lam,\,\mu\in P_+$.  
Inserting $z=1$ into \rf[e-hlam] we get $F_{\lam\mu}=e_{\lam\mu}(q,1)$.
Since $e_{\lam\mu}(t,1)$ is a rational function of $t$ and $q$ is a
transcendental root of it, $e_{\lam\mu}(t,1)\equiv 0$. In particular
$0=\lim_{t\to1}(t-t^{-1})^{-1}e_{\lam\mu}(t,1)=|\mu|-|\lam|$. Set
$m:=|\lam|=|\mu|$. Further we have $0=\lim_{t\to
  1}(t-t^{-1})^{-2}e_{\lam\mu}(t,1)$.  Since
\begin{align*}\tag{$*$}
\begin{split}
(t-t^{-1})^{-2}e_{\lam\mu}(t,1)=&\sum_{x\in\mu}\sgn(x)t^{\frac{N}{2}+c(x)}
\qnum[{\ts\frac{N}{2}}+c(x)]_t-
\\
&-\sum_{x\in\lam}\sgn(x)t^{-\frac{N}{2}-c(x)}\qnum[{\ts -\frac{N}{2}}-c(x)]_t
\end{split}
\end{align*}
the limit $t\to1$ gives $0=mN+c(\mu)+c(\lam)$.
Let us define $\widetilde{E}_\nu^+:=(2N)^{-1}(|\nu|N^2+2c(\nu)N-|\nu|^2)$ for
$\nu\in P_+$. By \rf[e-lim2], 
$\widetilde{E}_\nu^+\ge 0$ for $\nu\in P_{++}^0$ and $\widetilde{E}_\nu^+=0$ if and only if
$\nu=(0)$. It is easy to check that $\widetilde{E}^+_{\nu+(1^N)}=\widetilde{E}_\nu^+$ for
$\nu\in P_+$. Hence $\widetilde{E}_\nu^+\ge0$, $\nu\in P_+$, and $\widetilde{E}_\nu^+=0$ if and only if  $\nu=(n^N)$
for some $n\in\Z$. We conclude  that $c(\lam),c(\mu)\ge(2N)^{-1}(m^2-mN^2)$.
Inserting this into equation ($*$) we obtain $0\ge\frac{1}{N}m^2\ge0$, 
where equality holds on the left hand side if and only if $\lam=(n^N)$ and
$\mu=(l^N)$ for some $n,\,l\in\Z$ and on the right hand side if and only if
$m=0$. Hence $\lam=\mu=(0)$ and $z=1$ is regular. 

Finally consider the case when $z^Nq^{-2}=\zeta $ and $\zeta $ is a primitive
$m^{\,\mathrm{th}}$ root of unity ($m\in\N$). Let
$\tilde{z}=(t^2\zeta )^{1/N}$.
Observe that $e^\tau_{\lam+(m^N)}(t,\tilde{z})=e^\tau_{\lam}(t,\tilde{z})$
by \rf[e-recursion] and since $\zeta ^m=1$. Hence
$E_{\lam+(m^N),\mu+(m^N)}=F_{\lam+(m^N),\mu+(m^N)}=F_{\lam\mu}$
for $\lam ,\mu \in P_+$. Therefore the if part of the assertion
holds. Moreover, for the only if part we can assume that
$\lam,\mu \in P_{++}$.
Since $q$ is a transcendental root of the algebraic
function $f(t):=e_{\lam\mu}(t,\tilde{z})$, we conclude that $f(t)$ has to be
identically zero.
Therefore $\lim _{t\to \pm 1}f(t)$ has to be zero,
and by \rf[e-elampm] it follows that $|\mu|,|\lam|\in mN\Z$.
Then $\lim _{t\to 1}f(t)/(t^2-1)=0$. Further, one can compute
that $\lim _{t\to 1}t^{N-1}f(t)/(t^2-1)^2=\widetilde{E}_\mu^+
+\widetilde{E}_\lam^+$ (see \rf[e-lim2]). This sum is positive
except for the case $\widetilde{E}_\mu^+=\widetilde{E}_\lam^+=0$.
Moreover, the latter is equivalent to $\lam =(k^N),\mu =(l^N)$ for some
$k,l\in \Z$. Together with $|\mu|,|\lam|\in mN\Z$ we get the assertion.

(iii) There
exists an isomorphism  of Hopf bimodules $\gpz$ and $\gmz$, see
\cite[Subsection~14.6.1]{b-KS}.
In particular $X^+_0=X^-_0$ and consequently $E^+_\lam=E^-_\lam=:E_\lam$. 
By  \rf[e-X] and the definition of $\bo \ell^\pm$
one obtains
\begin{multline*}
\sss X_0^+(u^{i_1}_{k_1}\cdots u^{i_m}_{k_m}P_\lam{}^{k_1\cdots k_m}_{j_1\dots
  j_m})=
\\
=z^m(\Rda{m}\Rda{m-1}  \cdots\Rda{2}  \Rda[\,2]{1}\Rda{2}\cdots\Rda{m}-I)^{i_1\cdots
  i_m\, i}_{k_1\cdots k_m\, k}(D^{-1})^k_i \,P_\lam{}^{k_1\cdots
  k_m}_{j_1\cdots j_m},
\end{multline*}
where $P_\lam\in\Mor(\bo u^{\ot m})$. 
Paying attention to the $2$-admissible parameter $z\in\{-1,1\}$, the choice of the
coinvariant 1-form $\om_0^\tau$, and the definition of the $\sig$-metric it follows
from the remark after Proposition~7.2 in \cite{a-Heck1} that
$E_\lam=e_\lam(q)$, where
\begin{align*}
e_\lam(t):=\eps z^{|\lam|}(t-t^{-1})^2\sum_{x\in\lam}\qnum[N-\eps+2c(x)]_t.
\end{align*}
Suppose that $E_{\lam\mu}=0$. Since $\glm[\wedge]{\lam,\mu}=\{0\}$ for
$|\lam|\not\equiv|\mu|\mod (2)$ we may assume $|\lam|\equiv|\mu|\mod (2)$. Since $e_\lam(t)+e_\mu(t)$ is a rational function with
transcendental root $t=q$, $e_\lam+e_\mu\equiv 0$. In particular
\begin{equation}\label{e-lim}
\begin{aligned}
0\overset{!}{=}&\lim_{t\to1}\eps
z^{|\lam|}(t-t^{-1})^{-2}(e_\lam(t)+e_\mu(t))\\
=&(|\lam|+|\mu|)(N-\eps)+2c(\lam)+2c(\mu).
\end{aligned}
\end{equation}
Using the inequality $2c(\nu)N\ge |\nu|^2-|\nu|N^2$, $\nu\in P_{++}$, from the
proof of (ii) it is easily seen that for $\lam\ne(0)$ in the case $\eps=-1$
and in the case ($\eps=1$ and
$|\lam|>N$) we have $2c(\lam)>-|\lam|(N-\eps)$. Thus, by \rf[e-lim],
$\lam=\mu=(0)$ or ($\eps=1$ and $|\lam|,|\mu|\le N$). In the latter case
$2c(\lam)\ge -|\lam|^2+|\lam|$ where equality holds if and only if $\lam=(1^k)$, $1\le
k\le N$. Inserting this into \rf[e-lim] yields
$\lam,\mu\in\{(0),(1^N)\}$. Indeed we have $E_{\lam\mu}=0$ in these cases.
\end{proof}

\subsection*{Proof of Theorem\,\ref{t-co}} 
Let $[\rho ]$ denote the  cohomology class of  $\rho\in\gd$. Since $\dd$ is
bicovariant and $h^\lam\ast a=a\ast h^\lam$, the differential $\dd$ commutes
with $h^{\lam\mu}$, $\lam,\mu\in\AD$. In particular $h^{0,0}$
factorises to $h^{0,0}{}^\ast\colon \deR (\gd)\to \deR (\gdi)$. Since
$h^{0,0}\iota=\id$, where $\iota$ is the embedding of $\gdi$ into $\gd$,
$h^{0,0}{}^\ast$ is surjective. We prove injectivity. Let $\rho\in
\deR (\gd)$, in particular $\dd \rho=0$,
and suppose that $h^{0,0}{}^\ast(\rho)=[\rho^{0,0}]=0$.

(i) Consider first the case (a) of Theorem\,\ref{t-co} with $\Gamm=\gpz$. By
Lemma\,\ref{l-nonzero}, $E_{\lam \mu}\ne0$ if and only
if $(\lam,\mu)\ne(0,0)$. {}From \rf[e-laplace] and \rf[e-lrho] we obtain
\begin{align}
\begin{split}
\rho^{\lam\mu}&=(-1)^kE_{\lam\mu}^{-1}(-\dd\,\kod[+]+\kod[+]\,\dd)\,\rho^{\lam\mu},
\\
\rho^{\lam\mu}&=(-1)^kE_{\mu\lam}^{-1}(-\dd\,\kod[-]+\kod[-]\,\dd)\,\rho^{\lam\mu}
\end{split}
\label{e-hodge1}
\end{align}
for $(\lam,\mu)\ne(0,0)$. Since  $\dd$ commutes
with $h^{\lam\mu}$  we get
$\dd\rho^{\lam\mu}=h^{\lam\mu}(\dd \rho)=h^{\lam\mu}(0)=0$.
By \rf[e-hodge1], 
 $[\rho^{\lam\mu}]=(-1)^kE_{\lam\mu}^{-1}[-\dd\,\kod[+]\rho^{\lam\mu}]=0$  for
 $(\lam,\mu)\ne(0,0)$ (coboundary). Hence $[\rho]=[\rho^{0,0}]+\ds\sum_{(\lam,\mu)\ne(0,0)}[\rho^{\lam\mu}]=0$ and
$h^{0,0}{}^\ast$ is injective. In the same way the restrictions of the map
$h^{0,0}$ to the subspaces $\gdl$ and $\gdr$ yield isomorphisms
$\deR (\gdl)\cong
\deR (\gdi)$ and $\deR (\gdr)\cong \deR (\gdi)$, respectively. This proves (a)
in case  $\Gamm =\gpz$. For
$\Gamm =\gmz$ use
$\rho^{\lam\mu}=(-1)^kE_{\lam\mu}^{-1}(-\dd\,\kod[-]+\kod[-]\,\dd)\rho^{\lam\mu}$
and
$\rho^{\lam\mu}=(-1)^kE_{\mu\lam}^{-1}(-\dd\,\kod[+]+\kod[+]\,\dd)\rho^{\lam\mu}$
to get the same result.

(ii) Observe that for $a\in\A$ with $\dd a=0$ we have 
\begin{equation}\label{e-al}
\deR (a\Lamb)\cong a \deR (\Lamb).
\end{equation}
Consider the quantum group $\glqn$ and suppose that the parameter $z$ satisfies
the condition $z^Nq^{-2}=\zeta $, where $\zeta $ is a primitive
$m^{\,\mathrm{th}}$ root of unity, $m\in \N $.
Note that $\dd (\dett ^m)=0$ by \rf[e-dd]. Further we have
$\gl=\Gamm(0,0)\oplus\Gamm(0,(1,0,\dots,0,-1))$ by \rf[e-rightcoaction]
and $\bo u\ot\bo u^\fettc\cong\bo 1\oplus \bo u^{(1,0,\dots,0,-1)}$. By the
Littlewood-Richardson rule for tensor product representations of
$\mathrm{GL}(N)$ a necessary condition for $\glm[\wedge]{0,\mu}\ne\{0\}$ is
$|\mu|=0$. Since $\gd=\A \,\gdl$, $\glm[\wedge]{\lam,\mu}\ne\{0\}$ implies
$|\lam|=|\mu|$. Combining this with Lemma\,\ref{l-nonzero},  $E_{\lam\mu}=0$
is an eigenvalue of $\Lappm $ if and only if $\lam=\mu=(n^N)$ for some
$n\in m\Z$. Similarly as in (i) it follows that 
\[
h_\dett:=\sum_{n\in m\Z}h^{(n^N),(n^N)}\colon\gd\to
\bigoplus_{n\in m\Z}\glm[\wedge]{(n^N),(n^N)}
=\bigoplus_{n\in m\Z}\dett^{n}\gdi
\]
defines an isomorphism $h_\dett^\ast\colon
\deR (\gd)\to\bigoplus_{n\in m\Z}\deR (\dett^n \gdi)$. By \rf[e-al],
$\deR (\dett^{n} \gdi)=\dett^{n}\deR (\gdi)$ for $n\in m\Z$.
Since the images of both mappings $h_\dett\uhr\gdl$ and $h_\dett\uhr\gdr$
belong to $\gdi$, they define quasi-isomorphisms from $\gdl$ to $\gdi$ and from
$\gdr$ to $\gdi$, respectively. This proves the $\glqn$ part of (b).

(iii) Consider now the $\oqn$ case. Since $\gpz\cong\gmz$, it suffices to
carry out the proof for the calculus $\Gamm =\Gamm_{+,z}$, $z\in\{-1,1\}$. 
By Lemma\,\ref{l-nonzero}, $E^{\lam\mu}=0$ is an eigenvalue of
$\Lapp[]=Lapm[]$ if
and only if $\lam,\mu\in\{(0),(1^N)\}$. Similarly as in (i)
$h_\dett:=h^{0,0}+h^{(0),(1^N)}+h^{(1^N),(0)}+h^{(1^N),(1^N)}$ defines a
quasi-isomorphism of $\gd$ onto $\glm[\wedge]{\bo 1,\bo 1}+\glm[\wedge]{\bo
  1,\dett}+\glm[\wedge]{\dett,\bo 1}+\glm[\wedge]{\dett,\dett}$.
By the definition of $\bo \ell^\pm$, $\Rda{}$,
and $\dett$ we obtain $\ell^+{}^i_j(\dett)=x^{-N}\delta_{ij}$ and
$\ell^-{}^i_j(\dett)=y^N\delta_{ij}$. Consequently,
$X^+_{ij}(\dett)=(z^{-N}-1)\delta_{ij}$ and $\dd \dett=(z^{-N}-1)\dett
\om_0^+$. Hence $\dd\dett\ne 0$ in case of $\oqn[2n+1]$ and
$\Gamm_{\tau,-1}$. Otherwise $\dd\dett=0$.  
Since $\bo u\ot \bo
u^\fettc\cong\bo u^{(2)}\oplus \bo u^{(11)}\oplus\bo 1$, each irreducible subcorepresentation
of any tensor power $(\bo u\ott \bo u^\fettc)^{\ot k}$ corresponds to a Young
diagram  with  an {\em even}
number of boxes. Consequently, $\glm[\wedge]{\lam,\mu}=\{0\}$ if $|\lam|+|\mu|$
is odd.
In particular, $\glm[\wedge]{\bo 1,\dett}=\glm[\wedge]{\dett,\bo 1}=\{0\}$ if
$N$ is odd. For even $N$ however these spaces may be
nonzero. This completes the proof of Theorem\,\ref{t-co}.

\section{Proof of Theorem\,\ref{t-hodge}}\label{s-proof}

First we will show that the duality of $\dd_\tau$ and
$\kodpm_{-\tau}$ holds in a rather general setup.
Secondly, we will prove that for the quantum groups $\slqn$ and $\glqn$ the
differential calculi $\gdp$ and $\gdm$ are weakly isomorphic. Combining both
we obtain the proof of the second theorem.

\subsection*{Duality of differential and codifferential} 
\begin{proposition}\label{p-nondeg}
Suppose  that $\A$ is a cosemisimple Hopf algebra, $(\gp,\gm)$ is a dual
pair of bicovariant differential calculi and $\lam,\,\mu\in\AD$.
For $\nu \in \AD $ let
$\nu^\fettc\in\AD$ denote the class of the contragredient
corepresentation $(\bo u^\nu )^\fettc$.
\\
{\em (i)} The map
$h{\kri}\ctrpm{\cdot}{\cdot}\colon\glmt{\lam,\mu}\times\glmmt{\lam^\fettc,\mu^\fettc}\to\C$,
$\tau\in\{+,-\}$, $k\ge0$, is  non-degenerate. 
\\
{\em  (ii)} The restricted differential
$\dd_\tau\colon\glmt{\lam,\mu}\to\glmt[k+1]{\lam,\mu}$
is the dual operator to the restricted codifferential
$\kodpm_{-\tau}\colon\glmmt[k+1]{\lam^\fettc,\mu^\fettc}\to
\glmmt{\lam^\fettc,\mu^\fettc}$
with respect to the pairing $h{\kri}\ctrpm{\cdot}{\cdot}$.
\end{proposition} 
\begin{proof} (i) Set $\gd:=\gdt$ and fix $\lam\in\AD$, $k\in\N$. Let us prove
that $\sum_{\mu\in\AD}\glm{\lam,\mu}$ is finite
dimensional. The space $\gdl[k]\equiv\sum_{\nu\in\AD}\glm{0,\nu}$ is finite
dimensional since $\gl$ is. Suppose that
$\rho\in\sum_{\mu\in\AD}\glm{\lam,\mu}$. Since
$\rho=\rho_{(-2)}\cdot S(\rho_{(-1)})\rho_{(0)}$, we deduce that
$\rho\in\CC(\bo u^\lam)\gdl[k]$.
Because $\CC(\bo u^\lam)$ is finite dimensional,
the assertion follows.
Similarly, $\dim\sum_{\lam\in\AD}\glm{\lam,\mu}<\infty$.
For $\nu,\,\kappa\in\AD$ let $\{\rho^i_j\}$ and
$\{\zeta^m_n\}$ denote the linear bases of $\glmt{\lam,\mu}$ and
$\glmmt{\nu,\kappa}$,  respectively. Then we have 
\begin{align*}
\dlra(\rho^i_j)&=u^\lam_{ix}\ot\rho^x_y\ot u^\mu_{yj},
\\
\dlra(\zeta^k_l)&=u^\nu_{ka}\ot \zeta^a_b\ot u^\kappa_{bl}.
\end{align*}
Set $h^{ik}_{jl}:=h\ctrpm{\rho^i_j}{\zeta^k_l}$. By the left covariance of
$\ctrpm{\cdot}{\cdot}$ it follows that
\begin{equation*}
\begin{split}
u^\lam_{ix}u^\nu_{ky}h\ctrpm{\rho^x_j}{\zeta^y_l}&=(\id\ot
h)\Delta\ctrpm{\rho^i_j}{\zeta^ k_l}=1{\cdot}h\ctrpm{\rho^i_j}{\zeta^k_l},
\\
u^\lam_{ix}u^\nu_{ky}h^{xy}_{jl}&=h^{ik}_{jl}1.
\end{split}
\end{equation*}
Hence $(h^{ik}_{jl})_{i,k}\in\Mor(\bo 1,\bo u^\lam\ot\bo u^\nu)$ for all
$j,l$. 
By Schur's lemma we obtain $(h^{ik}_{jl})_{i,k}=0$ for
$\bo u^\nu\not\cong \bo u^{\lam^\fettc}$. Using right covariance, in a similar
way we get $(h^{ik}_{jl})_{j,l}\in\Mor(\bo u^\mu\ot\bo u^\kappa,\bo 1)$
for all $i,k$. Again by Schur's lemma $(h^{ik}_{jl})_{j,l}=0$ for $\bo
u^\kappa\not\cong \bo u^{\mu^\fettc}$. Suppose now that
$h\ctrpm{\rho}{\zeta}=0$ for a fixed $\rho\in\glmt{\lam,\mu}$ and all
$\zeta\in\glmmt{\lam^\fettc,\mu^\fettc}$. By the above arguments
$h\ctrpm{\rho}{\zeta}=0$ for all $\zeta\in\glmmt{\nu,\kappa}$,
$\nu,\kappa\in\AD$, i.\,e. for all $\zeta\in\Gaw[k]{-\tau}$. 
Since the Haar
functional is regular, i.\,e. $h(ab)=0$ for all $a\in\A$ implies $b=0$ and
$h(ab)=0$ for all $b\in\A$ implies $a=0$, and
since the pairing
$\ctrpm{\cdot}{\cdot}\colon\gdt[k]\ota\gdt[k]\to \A$ is non-degenerate, the pairing
$h\kri \ctrpm{\cdot}{\cdot}\colon\Gaw[k]{\tau}\ota\Gaw[k]{-\tau} \to \C$ is
also non-degenerate, cf. \cite[Section~6]{a-Heck1}.
Therefore $\rho=0$. Non-degeneracy in
the second component can be proved similarly. 

(ii) Suppose that  $\rho\in\glmt{\lam,\mu}$ and
$\zeta\in\glmmt[k+1]{\lam^\fettc,\mu^\fettc}$.  Because of  \rf[e-ctrlr], \rf[e-haar], the \mbox{$\sig$-symmetry}
of $g$, and since $\ctrpm{\rho}{\zeta}\in\Gamm_{-\tau}$, we obtain
\begin{align*}
h\ctrpm{\dd
  \rho}{\zeta}&=h\ctrpm{\om^\tau_0\land\rho-(-1)^k\rho\land\om^\tau_0}{\zeta}
\\
&=h\ctrpm{\om^\tau_0}{\ctrpm{\rho}{\zeta}}+(-1)^{k+1}h\ctrpm{\rho}{\ctrpm{\om^\tau_0}{\zeta}}
\\
&=h\ctrpm{\ctrpm{\rho}{\zeta}}{\om^\tau_0}+(-1)^{k+1}h\ctrpm{\rho}{\ctrpm{\om^\tau_0}{\zeta}}
\\
&=h\ctrpm{\rho}{\ctrpm{\zeta}{\om^\tau_0}}+(-1)^{k+1}h\ctrpm{\rho}{\ctrpm{\om^\tau_0}{\zeta}}
\\
&=h\ctrpm{\rho}{\kodpm \zeta}.
\end{align*}
The proof is complete.
\end{proof}
\subsection*{Homomorphy of differential calculi}\label{homo}
In this subsection we define and study the notion of homomorphic differential
calculi. 
Our aim is to show that for the quantum groups $\glqn$ and $\slqn$
the differential calculi $\gdp$ and $\gdm$ are weakly
isomorphic in the following sense. There exists a Hopf algebra automorphism
$F$ of
$\A$ which can be extended to a graded algebra isomorphism
$F\colon\gdp\to\gdm$ such that  $F\dd_+=\dd_- F$.

Suppose that $f\colon \BB\to 
\CC$ is an algebra homomorphism and $\Lam$ is a $\CC$-bimodule.
Then $\Lam$ is a $\BB$-bimodule via $a\cdot\eta\cdot b:=f(a)\eta f(b)$.
\begin{definition}
(i) Let $\BB$ and $\CC$ be  algebras and let $(\Gamm,\dd_1)$ and
$(\Lam,\dd_2)$ be first order differential calculi over $\BB$ and $\CC$,
respectively. The pair $(f,F)$ of an algebra
homomorphism $f\colon\BB\to \CC$ and a $\BB$-bimodule homomorphism
$F\colon\Gamm\to\Lam$ is called a {\em homomorphism of the first order
differential calculi $\Gamm$ and $\Lam$} if 
\begin{align}\label{e-hom1}
F\dd_1=\dd_2 f.
\end{align}

(ii) Let $\gd$ and $\ld$ be differential calculi over $\BB$ and $\CC$,
respectively. A
graded algebra homomorphism $F\colon\gd\to\ld$ is called a {\em homomorphism of
the   differential calculi} $\gd$ and $\ld$ if
\begin{align}\label{e-hom2}
F\dd_1=\dd_2 F.
\end{align}
\end{definition}
Recall  that two  differential calculi $\gd$ and $\ld$ over the same  algebra $\BB$  are isomorphic {\em in the strong
sense} if and only if there exists a bijective homomorphism $F\colon\gd\to\ld$
of differential calculi with  
$F_0=\id$.

The next lemma characterises homomorphic differential calculi in terms of
their associated right ideals and in terms of their quantum tangent spaces.
\begin{lemma}\label{l-hom}
Let $\Gamm$ and $\Lam$ be left-covariant first order differential calculi over
the Hopf algebras $\A$ and $\BB$ with associated right ideals $\RR$ and $\SS$
and quantum tangent spaces $\XX$ and $\YY$, respectively. Suppose that $f\colon\A\to
\BB$ is a Hopf algebra homomorphism. The following are equivalent:
\\
{\em (i)} $f(\RR)\subseteq\SS.$
\\
{\em (ii)} $f^\fettt(\YY)\subseteq\XX.$
\\
{\em (iii)} There exists a unique homomorphism $(f,F)$ of
 the  left-covariant FODC $\Gamm$ and $\Lam$.
\end{lemma}
\begin{proof}
(iii)$\to$(i). Fix $r\in\RR$. Since $f$ is a Hopf algebra homomorphism and
$F(a\dd b)=f(a)\dd f(b)$, we have
\begin{equation}
\begin{aligned}
\om_\Lam(f(r))&=S\bigl(f(r)_{(1)}\bigr)\dd f(r)_{(2)}
=S(f(r_{(1)}))\dd f(r_{(2)})
\\
&=f(Sr_{(1)})F(\dd r_{(2)})=F(\om(r))=0.
\end{aligned}
\end{equation}
Moreover, $\ve(f(r))=\ve(r)=0$. Hence $f(r)\in\SS$.

(i)$\to$(ii). Since $\SS$ and $\YY$ are orthogonal subspaces
with respect to the pairing of $\A $ and $\A^\circ $,
we have $f^\fettt(Y)(r)=Y(f(r))\in Y(\SS)=\{0\}$
for all $r\in\RR$ and $Y\in\YY$.
Furthermore, $f^\fettt(Y)(1)=Y(f(1))=Y(1)=0$. Hence $f^\fettt(Y)\in\XX$.

(ii)$\to$(iii).   $F$ is uniquely determined by \rf[e-hom1], since
$F(a_i\dd b_i):=f(a_i)\dd f(b_i)$. We show that $F$ is well-defined. Let
$\{\om_i\}$ and $\{\eta_j\}$ be linear bases of $\gl$ and $\Lam_\ll$, 
respectively, and let $\{X_i\}$ and $\{Y_j\}$ be the 
corresponding dual bases of $\XX$ and $\YY$, respectively. By assumption there 
exist $\alpha_{ij}\in\C$ such that
$f^\fettt(Y_j)=\alpha_{ij}X_i$. Suppose that $a_i\dd
b_i=0$. Then we have $0=a_i(X_k\ast b_i)\om_k$ and consequently $0=a_i(X_k\ast
b_i)$ for all $k$. Using this fact, we conclude that
\begin{align*}
f(a_i)\dd
f(b_i)&=f(a_i)f(b_{i(1)})\om_\Lam(f(b_{i(2)}))=f(a_i)f(b_{i(1)})Y_j(f(b_{i(2)}))\eta_j
\\
&=f(a_ib_{i(1)})f^\fettt(Y_j)(b_{i(2)})\eta_j=\alpha_{kj}f(a_ib_{i(1)})
X_k(b_{i(2)})\eta_j
\\
&=\alpha_{kj}f(a_i(X_k\ast b_i))\eta_j=0.
\end{align*}
Hence  $F$ is well-defined. 
\end{proof}
The next lemma is straightforward to prove using covariance of $\dd$ and the
properties of $F_0$. We omit the proof.
\begin{lemma}\label{l-comod} Suppose that $\gd$ and $\ld$ are left-covariant
  differential calculi over the Hopf algebras $\A$ and $\BB$, respectively. 
Let $F\colon\gd\to\ld$ be a homomorphism of  differential
calculi and $F_0$ a Hopf algebra homomorphism. Then we have $(F_0\ot
F)\dl=\dl F$ and $F(\rho\rac a)=F\rho\rac F_0a$ for $\rho\in\gd$ and
$a\in\A$.
\end{lemma}
Replacing left-covariance by right-covariance in the above lemma the first
assertion reads as $(F\ot F_0)\dr=\dr F$.
Now we shall apply the new notion to our main example.
\begin{proposition}\label{p-auto}
Let $G_q$ be one of the quantum groups $\glqn$ or  $\slqn$ and
$\A =\OO(G_q)$.
For $k=1,\dots,N$, set $k'=N+1-k$. 

{\em (i)} There exists a unique bijective homomorphism $F\colon\gdp\to\gdm$ of 
differential calculi such
that
\begin{align}\label{e-auto}
F(u^a_b)=Su^{b'}_{a'},\quad a,b=1,\dots,N.
\end{align}
{\em (ii)} For all $\lam,\mu \in \AD $ the restriction of $F$
to $\glmp{\lam,\mu}$ is a bijection onto $\glmm{\lam^\fettc,\mu^\fettc}$.
\end{proposition}
\begin{proof} (a) First it is to show that there exists a Hopf algebra
automorphism $F\colon\A\to\A$
which satisfies  \rf[e-auto]. To do this we prove that $F$ preserves the
relations of the Hopf algebra $\A$.
It is easily  shown that  $\Rda[ab]{rs}=\Rda[s'r']{b'a'}$ and
$d_k^{-1}=d_{k'}$. Moreover,  the $q$-antisymmetric tensor satisfies
$\ve_{i_1\cdots\, i_N}=\ve_{i_N^\prime\cdots\, i_1^\prime}$. 
We show that the algebra homomorphism
$f\colon \C\langle u^i_j\rangle\to \C\langle v^i_j\rangle$, given by
$f(u^a_b)=v^{b'}_{a'}$ maps the generating relations appearing in the
definition of the Hopf algebra $\A $ to those of the Hopf algebra
$\A^{\op,\cop}$.
Here we assume that $\C\langle v^i_j\rangle$ has both opposite multiplication
and opposite comultiplication.
By the above identity  for the matrix $\Rda{}$ we have
\[
f(\Rda[ab]{xy}u^x_r\ot u^y_s-u^a_x\ot u^b_y
\Rda[xy]{rs})=\Rda[y'x']{b'a'}v^{r'}_{x'}\ot v^{s'}_{y'}-
v^{x'}_{a'}\ot v^{y'}_{b'} \Rda[s'r']{y'x'}.
\]
The right hand side  generates the relations of the bialgebra $\A^{\op,\cop}$.
Similarly
one shows consistency with the $q$-determinant relation. Finally we have
$(f\ot f)\Delta u^a_b=v^{x'}_{a'}\ot v^{b'}_{x'}=\Delta
v^{b'}_{a'}=\Delta(f(u^a_b))$
and   $\ve(f(u^a_b))=\delta_{ab}=\ve(u^a_b)$. Hence $f$ is a homomorphism
of bialgebras. Since both $\A$ and $\A^{\op,\cop}$ are Hopf algebras,
$f$ is a Hopf algebra homomorphism. Since the antipode is a Hopf algebra map
of $\A^{\op,\cop}\to \A$, $F=S\kri f$ is a Hopf algebra automorphism. Its
inverse $F^{-1}$ is given by $F^{-1}(u^a_b)=S^{-1}(u^{b'}_{a'})$.  

Next we show that $F(\CC(\bo u^\lam))=\CC(\bo u^{\lam^\fettc})$. Let
$P^\lam\in\Mor(\bo u^{\ot k})$ be a primitive idempotent such that $\CC(\bo u^\lam)=\langle (P^\lam)^\nvec_\xvec\,
u^\xvec_\mvec\,|\,\nvec,\mvec\in\{1,\dots,N\}^k\rangle$. For
$\xvec=(x_1,\dots, x_k)$ we write
$\xvec\,^\prime=(x_1^\prime,\dots,x_k^\prime)$ and
$\xinvp=(x_k^\prime,\dots,x_1^\prime)$. Let us show that $(Q^\lam)^\nvec_\mvec:=(P^\lam)^\minvp_\ninvp$ is a projector equivalent to $P^\lam$. Idempotents $P$ and $Q$ are called equivalent, if
there exist $A,\,B\in\Mor(\bo u^{\ot k})$ such that $AB=P$ and $BA=Q$. For
this let $\alpha$ and $\beta$  denote the algebra automorphism and algebra antiautomorphism of $\Mor(\bo u^{\ot   k})$ defined by 
\begin{align*}
\alpha(\Rda{n,n+1})=\Rda{k-n,k-n+1}\quad \text{and}\quad\beta(\Rda{n,n+1})=\Rda{n,n+1},
\end{align*}
$n=1,\dots, k-1$, respectively.  By  the theory of Hecke algebras  it is easy to see that $\alpha$ and $\beta$ map
each twosided ideal of $\Mor(\bo u^{\ot k})$ into itself. In particular, the
image of a primitive idempotent is an equivalent primitive idempotent.
By induction on $k$ we will show that 
\begin{equation}
T\,^\minvp_\ninvp=\beta(\alpha(T))^\nvec_\mvec.
\end{equation}
It is well known that for   $T\in\Mor(\bo u^{\ot k})$ there exist
$X,\,Y\in\Mor(\bo u^{\ot 2})$ and $B\in\Mor(\bo u^{\ot k-1})$ such that
$T=X_{k-1,k}BY_{k-1,k}$. Since $X$ and $Y$ can be
written in terms of $\Rda{}$, $\Rda[-1]{}$, and $I$, and since 
$\Rda[ab]{rs}=\Rda[s'r']{b'a'}$, we have by induction assumption
\begin{align*}
T^{m_k'\cdots m_1'}_{n_k'\cdots n_1'}&=X^{m_2' m_1'}_{x'y'}B^{m_k'\cdots
  m_3'x'}_{n_k'\cdots n_3'z'}Y^{z'y'}_{n_2'n_1'}
=Y^{n_1n_2}_{yz}\beta(\alpha( B))^{z n_3\cdots n_k}_{x m_3\cdots
  m_k}X^{yx}_{m_1m_2}
\\
&=(Y_{12}\beta(\alpha(B)) X_{12})^{\nvec}_\mvec= \beta(\alpha(X_{k-1,k}BY_{k-1,k}))^\nvec_\mvec.
\end{align*}
 The character of $\bo u^\lam$ is
$\chi(\bo u^\lam)=(P^\lam)^\nvec_\xvec\, u^\xvec_\nvec$ \cite[Lemma~5.1]{a-Schueler1}. Since
$Q^\lam=\beta(\alpha(P^\lam))$ is an equivalent idempotent, the corresponding
characters of the corepresentations $(P^\lam)^\nvec_\xvec\, u^\xvec_\mvec$ and 
$(Q^\lam)^\nvec_\xvec\, u^\xvec_\mvec$ coincide \cite[Lemma~5.1]{a-Schueler1}.
We conclude that
\begin{align}
\chi(F\bo u^\lam)&=(P^\lam)^\nvec_\xvec\,(u^\fettc)^{\xvecp}_{\nvecp}
=(Q^\lam)^{\xinvp}_{\ninvp}\,S(u^{\ninvp}_{\xinvp})=S(\chi(\bo u^\lam))
=\chi(\bo u^{\lam^\fettc}).
\end{align}
Since $F$ is a Hopf algebra homomorphism, matrix elements of irreducible
corepresentations are mapped into each other. Hence $F(\CC(\bo
u^\lam))=\CC(\bo u^{\lam^\fettc})$.

(b) We show that $F^\fettt(\XX ^-)\subseteq \XX ^+$. Since $F$ is a coalgebra
homomorphism, $F^\fettt$ is multiplicative on the subalgebra of $\A^\circ$
generated by the  matrix elements $\ell^\pm{}^a_b$, $a,b=1,\dots,N$. We
compute $F^\fettt(\ell^\pm{}^a_b)$ on the generators of $\A$.
\begin{align*}
F^\fettt(\ell^\pm{}^a_b)(u^r_s)&=\ell^\pm{}^a_b(Fu^r_s)=\ell^\pm{}^a_b(Su^{s'}_{r'})
\\&=p^{\pm 1}(\Rda[\mp1]{})^{s'a}_{br'}=p^{\pm1}(\Rda[\mp1]{})^{rb'}_{a's}=\ell^\pm{}^{b'}_{a'}(Su^r_s)=(S\ell^\pm{}^{b'}_{a'})(u^r_s),
\end{align*}
where $p=x$ in the $\bo \ell^+$-case and $p=y$ in the $\bo \ell^-$-case.
Since both $F^\fettt(\ell^\pm{}^a_b)$ and  $S\ell^\pm{}^{b'}_{a'}$ are
representations of $\A$, we obtain
$F^\fettt(\ell^\pm{}^a_b)=S\ell^\pm{}^{b'}_{a'}$.
By \cite[Theorem~9.\,1]{a-FRT} we have
\begin{align*}
\ell^+{}^r_v\,S\ell^-{}^w_s\Rda[vb]{wc}=\Rda[rv]{sw}S\ell^-{}^b_v\,\ell^+{}^w_c.
\end{align*}
Multiplying this relation by $D^s_r$ and noting that $D^s_r\Rda[rv]{sw}=\rt
\delta_{vw}$, we get
\begin{align*}
D^s_r\ell^+{}^r_vS\ell^-{}^w_s\Rda[vb]{wc}=\rt S\ell^-{}^b_v\ell^+{}^v_c.
\end{align*}
Multiplying the latter by $\Rdam[mc]{na}(D^{-1})^a_bD^n_k$ and using
the identity
$\Rda[vb]{wc}\Rdam[mc]{na}(D^{-1})^a_bD^n_k=\delta_{vn}\delta_{wm}$
we obtain
\begin{align}\label{e-lpslm}
D^s_r\ell^+{}^r_kS\ell^-{}^m_s =\rt 
S\ell^-{}^b_v\ell^+{}^v_c\Rdam[mc]{na}(D^{-1})^a_bD^n_k.
\end{align}
By  \rf[e-X], \rf[e-lpslm],  $\Rdam[jc]{na}(D^{-1})^a_c=\rt^{-1}\delta_{jn}$,
and $d_k^{-1}=d_{k'}$, we then have
\begin{align*}
F^\fettt(X ^-_{i'j'})&=
F^\fettt\bigl((D^{-1})^k_lS^{-1}(\ell^+{}^{i'}_k)\ell^-{}^l_{j'}\bigr)
-(D^{-1})^{i'}_{j'}
\\
&=(D^{-1})^k_lS^{-1}(S\ell^+{}^{k'}_{i})S\ell^-{}^j_{l'}-D^i_j
\\
&=D_{k'}^{l'}\ell^+{}^{k'}_iS\ell^-{}^j_{l'}-D^i_j
\\
&=\rt S\ell^-{}^b_v\ell^+{}^v_c\Rdam[jc]{na}(D^{-1})^a_bD^n_i-D^i_j
\\
&=\rt( X^+_{bc}+\delta_{bc})\Rdam[jc]{na}(D^{-1})^a_bD^n_i-D^i_j
\\
&=\rt X ^+_{bc}\Rdam[jc]{na}(D^{-1})^a_bD^n_i.
\end{align*}
This completes the proof of (b). Note that $F^\fettt\colon\XX ^-\to\XX ^+$ is
bijective since $X ^+_{ab}=\rt^{-1}
F^\fettt(X ^-_{i'j'})(D^{-1})^i_k\Rda[ka]{lb}D^l_j$. By Lemma\,\ref{l-hom},
$F(a\dd b):=F(a)\dd F(b)$ is a well-defined $\A$-module map from $\gp$ to
$\gm$. Similarly, $F^{-1}\colon\gm\to\gp$, $F^{-1}(a\dd b):=F^{-1}(a)\dd
F^{-1}(b)$ is a well-defined $\A$-module map inverse to $F$. 

(c) Consider now higher order forms.
Let $\tau \in \{+,-\}$ and $F^\tau $ denote
$F$ for $\tau =+$ and $F^{-1}$ for $\tau =-$.
Since $F^\tau$ is an $\A$-bimodule map, we can extend
$F^\tau $ to an algebra map
$F^\tau \colon\gt^\ot\to \gmt^\ot$.
We prove that $F\sig=\sig F$ in $\Gamm ^{\ot 2}_+$. Since both $F$ and
$\sig$ are $\A$-bimodule maps and $\gp \ota \gp $ is a free left
$\A$-module with basis $(\gp \ota \gp )_\ll$, it suffices to prove
this equation for left-coinvariant elements. Let $\rho,\,\xi\in(\gp)_\ll$. By
\rf[e-sigli] and Lemma\,\ref{l-comod} we have
\begin{align*}
F\sig(\rho\ota\xi)&=F(\xi_{(0)}\ota(\rho\rac\xi_{(1)}))=F(\xi_{(0)})\ota
F(\rho\rac\xi_{(1)})
\\
&=(F\xi)_{(0)}\ota (F\rho\rac(F\xi)_{(1)})=\sig(F\rho\,\ota F\xi).
\end{align*}
Hence $F$ commutes with the antisymmetriser $A_k$, $k\in\N$. The same is true
for $F^{-1}$. Consequently, $F^\tau\colon\gdt\to\Gamm^\land_{-\tau}$
is a well-defined algebra map and $F \,F^{-1}=F^{-1}F=\id$.
Now let us prove assertion (ii). Let
$\rho\in\glmp{\lam,\mu}$. By Lemma\,\ref{l-comod} and the last part
of (a), $F\rho\in\glmm[k]{\lam^\fettc,\mu^\fettc}$. Since
$F^{-1}\colon\gdm\to\gdp$ also satisfies the assumptions of
Lemma\,\ref{l-comod}, it follows immediately that $F^{-1}\uhr
\glmm[k]{\lam^\fettc,\mu^\fettc}$ is inverse to $F\uhr\glmp[k]{\lam,\mu}$.
\end{proof}

\begin{remark}
For the B-, C-, and D-series the differential calculi $\gdp$ and $\gdm$ are
isomorphic (in the strong sense: $F_0=\id$).
\end{remark}
Now we are able to finish the proof of Theorem\,\ref{t-hodge}.  By
\cite[Theorem~3.2\,(iii)]{a-Schueler1} the differential $\dd$ vanishes on
$\gdi[k]$. Hence $\deR ^k(\gdi)\cong \gdi[k]$. Combining this with
Theorem\,\ref{t-co} gives \rf[e-cohom].
Since $\gdi[k]=\glm[k]{0,0}$ it follows from
\rf[e-dec] that $\deR ^k(\gd)$ is a direct summand in \rf[e-hodge].
By \rf[e-hodge1], for $(\lam,\mu)\ne(0,0)$ we have the following formulae:
\begin{equation}
\begin{split}\label{e-hodge2}
\glmt{\lam,\mu}&=\dd\glmt[k-1]{\lam,\mu}+\partial^+\glmt[k+1]{\lam,\mu },
\\
\glmt{\lam,\mu}&=\dd\glmt[k-1]{\lam,\mu}+\partial^-\glmt[k+1]{\lam,\mu }.
\end{split}
\end{equation}
We have to prove that both sums are direct.  Since all vector spaces
appearing in \rf[e-hodge2] are finite dimensional,
it suffices to compare their dimensions. 
We denote the restriction of a linear map  $f\colon\gd\to\gd$ to the space
$\glm[k]{\lam,\mu}$ by $f^{k,\lam,\mu }$. By Proposition\,\ref{p-nondeg}
we have $\rank\kodpm_\tau{}^{,k+1,\lam,\mu}=\rank
\dd_{-\tau}^{k,\lam^\fettc,\mu^\fettc}$.
Indeed, both $\dd_{-\tau }^{k,\lam^\fettc,\mu^\fettc}$ and
$\kodpm_\tau{}^{,k+1,\lam,\mu}$ are linear mappings acting on finite
dimensional vector spaces and they are dual to each other.
By Proposition\,\ref{p-auto},
$\dd_{-}^{k,\lam^\fettc,\mu^\fettc} F=F \dd_+^{k,\lam,\mu} $
and $F$ is bijective. We conclude that
$\rank\kodpm_\tau{}^{,k+1,\lam,\mu}=
\rank\dd_{-\tau}^{k,\lam^\fettc,\mu^\fettc}=
\rank\dd_\tau^{k,\lam,\mu}$. Since  $\deR ^k(\gdt(\lam,\mu))=\{0\}$,
$\dim\ker\dd_\tau^{k,\lam,\mu}=\rank\dd_\tau^{k-1,\lam,\mu}$. Finally we obtain
\begin{align*}
\dim \glmt{\lam,\mu}&=\dim \ker\dd_\tau^{k,\lam,\mu}+\rank\dd_\tau^{k,\lam,\mu}
\\
&=\rank\dd_\tau^{k-1,\lam,\mu}+\rank\kodpm_\tau{}^{,\,k+1,\lam,\mu}.
\end{align*}
It follows that the sums \rf[e-hodge2] are direct. The proof of
Theorem\,\ref{t-hodge} is complete.
\subsection*{Acknowledgement} We are grateful to Konrad Schm\" udgen for
suggesting this problem and for helpful comments.
\bibliographystyle{axel}
\bibliography{bibdat}

\affil{Department of Mathematics\\
University of Leipzig\\
Augustusplatz 10\\
04109 Leipzig\\
Germany}

\email{heckenbe,schueler@mathematik.uni-leipzig.de}

\end{article}
\end{document}